\documentclass[12pt,reqno]{amsart}
\usepackage{tikz-cd}
\usepackage{amscd}
\usepackage{amsmath, amsfonts,amssymb,amsthm,mathrsfs,xypic}
\usepackage{graphicx}
\usepackage{fancybox}
\calclayout
\usepackage{amsmath}
\usepackage{hyperref}
\usepackage{cleveref}
\usepackage[all]{xy}
\usepackage{extarrows}

\usepackage{mathtools}
\usepackage{mathrsfs,stmaryrd,wasysym}
\usepackage{graphics}
\usepackage{xypic}
\usepackage{multirow}
\usepackage{makecell}

\usepackage[english]{babel}
\usepackage{amsmath}

\allowdisplaybreaks[4]

\usepackage{cases}

\usepackage{dsfont}
\usepackage{color}
\SelectTips{cm}{}
\calclayout

\numberwithin{equation}{section}

\theoremstyle{plain}
\newtheorem{theorem}[equation]{Theorem}
\newtheorem{proposition}[equation]{Proposition}
\newtheorem{lemma}[equation]{Lemma}
\newtheorem{corollary}[equation]{Corollary}

\theoremstyle{definition}
\newtheorem{definition}[equation]{Definition}

\theoremstyle{remark}
\newtheorem{remark}[equation]{Remark}

\newcommand{\xrto}{\xrightarrow}

\def\mcP{\mathcal{P}}
\def\mcE{\mathcal{E}}

\def\mcL{\mathcal{L}}

\def\nslash{\:\notslash\:}
\def\mcH{\mathcal H}
\def\mcC{\mathcal C}

\def\mfD{\mathfrak{D}}
\def\mcW{\mathcal W}
\def\parti{\partial}

\def\bi{\text{\boldmath$i$}}

\def\bj{\text{\boldmath$j$}}

\begin{document}
	
	\title[\tiny{The cyclotomic non-degenerate Hecke algebras}]{The cyclotomic non-degenerate Hecke algebras of arbitrary Weyl groups}
	\author [\tiny{Fan Kong and Zhi-Wei Li}] {Fan Kong and Zhi-Wei Li*}
	
	\begin{abstract}
We construct Khovanov--Lauda--Rouquier-like generators for certain modified forms of non-degenerate affine Hecke algebras associated with arbitrary Weyl groups after localization. These generators yield non-graded KLR-like presentations of the modified algebras. As an application, we establish isomorphisms between direct sums of blocks of cyclotomic Hecke algebras of arbitrary Weyl groups and cyclotomic quotients of the resulting non-graded KLR-like algebras. We also explain how these isomorphisms identify the generalized weight-space decomposition of cyclotomic Hecke modules with the idempotent decomposition of modules over the KLR-like presentation.
\end{abstract}
	
	\date{\today}
	\thanks{*Corresponding author}
	\thanks{The first author is supported by the NSF of China (No. 11501057); the second author is supported by the NSF of China (No. 11671174).}
	
	\address{Fan Kong\\ School of Mathematics and Statistics \\
		Southwest University \\ Chongqing 400715\\ P. R. China.}
	\email{kongfan85@126.com}

	\address{Zhi-Wei Li\\ School of Mathematics and Statistics \\
		Jiangsu Normal University\\	Xuzhou 221116 Jiangsu \\ P. R. China.}
	\email{zhiweili@jsnu.edu.cn}
	\maketitle
	
	\setcounter{tocdepth}{1}
	
	\section{Introduction}

This paper is the non-degenerate counterpart of our previous work on cyclotomic degenerate Hecke algebras of arbitrary Weyl groups \cite{Kong-Li}. Affine Hecke algebras may be viewed as deformations of group algebras of affine Weyl groups, and the Hecke algebra of the finite Weyl group occurs naturally as a subalgebra. They arise in the representation theory of reductive $p$-adic groups \cite{Lusztig89,Lusztig88} and are closely related to several areas of representation theory and mathematical physics; see, for example, \cite{Cherednik91,Cherednik05,Mac}.

Cyclotomic quotients provide an important tool for studying finite-dimensional representations of affine Hecke algebras. In type $A$, Brundan and Kleshchev \cite{Brundan-Kleshchev09} constructed explicit isomorphisms between blocks of cyclotomic Hecke algebras for symmetric groups and cyclotomic Khovanov--Lauda--Rouquier (KLR) algebras. Independently, Rouquier obtained a closely related result in \cite[Section 3.2]{Ro} by a different approach. KLR algebras, introduced by Khovanov and Lauda \cite{KL} and Rouquier \cite{Ro}, are graded algebras designed to categorify quantum groups. Their representation theory in type $A$ is closely connected with that of affine Hecke algebras of symmetric groups \cite{Ro}, and the Brundan--Kleshchev--Rouquier (BKR) construction provides a powerful way to pass between the two settings \cite{BKW}.

The main purpose of this paper is to construct non-graded KLR-like presentations associated with non-degenerate affine Hecke algebras of arbitrary Weyl groups. To do this, we introduce, after localization, certain modified forms of affine Hecke algebras obtained by adjoining a family of orthogonal idempotents indexed by a Weyl-group orbit. In types $B$ and $D$, these modified algebras are related in spirit to the graded algebras introduced by Varagnolo and Vasserot \cite{VV}, Shan--Varagnolo--Vasserot \cite{SVV}, and their generalizations in \cite{Andecy-Walker17,Andecy-Walker19}. A crucial difference is that the KLR-like algebras arising in the present paper are, in general, not $\mathbb{Z}$-graded.

The affine Hecke algebra used here differs from the algebra appearing in the type-$A$ Brundan--Kleshchev construction \cite{Brundan-Kleshchev09}, although in type $A$ it occurs naturally as a subalgebra of the latter. The relationship between the corresponding cyclotomic quotients is not immediate. Moreover, the cyclotomic degenerate and non-degenerate Hecke algebras in the usual Brundan--Kleshchev--Rouquier setting are related to the same cyclotomic KLR algebra, whereas the analogous algebras constructed in our degenerate and non-degenerate settings are different. Consequently, the type-$A$ specialization of the isomorphisms obtained here does not recover the Brundan--Kleshchev--Rouquier isomorphisms.

Intertwining elements are the main technical tool in our construction. Following Cherednik \cite{Cherednik91}, we work with intertwiners in a localization of the affine Hecke algebra and then modify them in a way inspired by the Brundan--Kleshchev construction \cite{Brundan-Kleshchev09}. This produces Khovanov--Lauda--Rouquier-like generators for the modified forms of affine Hecke algebras. We derive a complete set of relations for these generators and thereby obtain non-graded KLR-like presentations.

As an application, we introduce cyclotomic quotients of the resulting KLR-like algebras and prove isomorphisms between them and direct sums of blocks of cyclotomic Hecke algebras of arbitrary Weyl groups. Thus the paper provides a uniform non-graded analogue, for arbitrary Weyl groups, of the passage from cyclotomic Hecke algebras to KLR-type presentations in type $A$.

The main result was announced in our preprint \cite{Kong-Li21}, posted on arXiv in 2021. We do not intend to publish that preprint separately.

The paper is organized as follows. In Section~2, we recall affine Hecke algebras of arbitrary Weyl groups and their localizations, review the relevant intertwining elements, and introduce the modified forms together with modified intertwiners. In Section~3, we construct KLR-like generators for these modified algebras and establish a complete set of defining relations. In Section~4, we introduce the corresponding KLR-like algebras and their cyclotomic quotients, construct the resulting non-graded Brundan--Kleshchev--Rouquier-like isomorphisms for cyclotomic Hecke algebras of arbitrary Weyl groups, and explain the representation-theoretic interpretation of the new presentation.

\section{Modified forms of affine Hecke algebras and intertwiners}

	\subsection{The Weyl groups} 	Let $\Sigma\subseteq \mathbb{R}^n$ be a root system of type $A_n,B_n,\cdots,G_2$, and let $\mfD$ be the Dynkin diagram of $\Sigma$ with vertex set $[n]=\{1,2,\cdots,n\}$. 
We use $C=(c_{ij})_{n\times n}$ to denote the corresponding Cartan matrix of $\mfD$, where 
		\begin{tiny}\begin{align*}c_{ij}=\begin{cases}2 & \text{if $i=j$};\\
				0& \text{if $i\nslash j$};\\
				-1 & \text{if $i-\-j$, or $\xymatrix@C10pt{i\ar@{=}[r]|{\langle}\ar@{=}[r]&j},$ or $\xymatrix@C10pt{i\ar@3{-}[r]|{\langle}\ar@3{-}[r]&j}$};\\
				-2 & \text{if$ \xymatrix@C10pt{i\ar@{=}[r]|{\rangle}\ar@{=}[r]&j}$};\\
				-3 & \text{if $\xymatrix@C10pt{i\ar@3{-}[r]|{\rangle}\ar@3{-}[r]&j}$}.
	\end{cases}\end{align*}\end{tiny}
	Here we use $i \nslash j$ to indicate that $i$ is not connected to $j$ by an edge, $i -\- j$ indicates that $i$ is connected with $j$ by an edge, $\xymatrix@C10pt{i\ar@{=}[r]|{\rangle}\ar@{=}[r]&j}$ indicates that $i$ is connected with $j$ by a double edge and there is an arrow from $i$ to $j$, and $\xymatrix@C10pt{i\ar@3{-}[r]|{\rangle}\ar@3{-}[r]&j}$ indicates that $i$ is connected with $j$ by a triple edge and there is an arrow from $i$ to $j$. To guarantee the correctness of Lemma \ref{lem:deynilpotent}, we assume that there is a directed path from each vertex to vertex $1$ in the Dynkin diagram $\mfD$.
	 
	For $\alpha\in \Sigma$, let $s_\alpha$ be the orthogonal reflection in the hyperplane $\{v\in \mathbb{R}^n\ | \ (\alpha, v)=0\}$ with respect to the standard Euclidean inner product $(\ ,\ )$ on $\mathbb{R}^n$. 
Let $\{\alpha_i\ | \ i\in[n]\}\subseteq \Sigma$ be the simple roots for some fixed Weyl chamber and $\mcW$ the Weyl group generated by the reflections $\{s_i:=s_{\alpha_i}\ | \ i\in [n]\}$. Then the order of $s_is_j$ with $i\neq j$ is $m_{ij}$, where 
\begin{tiny}	\begin{align*}m_{ij}=\begin{cases}2& \text{if $i\nslash j$};\\
			3 & \text{if $i-\-j$};\\ 
		4 & \text{if$ \xymatrix@C10pt{i\ar@{=}[r]&j}$};\\
		6 & \text{if $\xymatrix@C10pt{i\ar@3{-}[r]&j}$}.
\end{cases}\end{align*}\end{tiny}

The length of $w\in \mcW$ (that is, the length of any reduced expression of $w$ with respect to $\{s_i, i\in [n]\}$) is denoted by $\ell(w)$. 

	\subsection{The Demazure operators} Following \cite{Cherednik91}, let $\Bbbk[X;F]$ be the polynomial algebra generated by $\{X_\alpha\ | \ \alpha\in \Sigma\}$ over a fixed field $\Bbbk$ with the relations $$F\colon X_{\alpha+\beta}=X_{\alpha}X_{\beta} \ \mbox{and}\ X_0=1$$ for $\alpha, \beta\in \Sigma, \alpha+\beta\in \Sigma\cup \{0\}$. Then there is an action of the Weyl group $\mcW$ on $\Bbbk[X;F]$ via 
	$$w(X_\alpha)=X_{w(\alpha)}, \ w\in \mcW, \alpha\in \Sigma.$$
	It extends to the field of fractions $\Bbbk(X;F)$ of $\Bbbk[X;F]$ by defining $$w(\tfrac{f}{g})=\tfrac{w(f)}{w(g)}$$ for $w\in \mcW$ and $\tfrac{f}{g}\in \Bbbk(X;F)$. 
	Set $X_i=X_{\alpha_i}$. Then
	$$s_i(X_{j})=X_{s_i(\alpha_j)}=X_{\alpha_j-c_{ji}\alpha_i}=X_{j}X_{i}^{-c_{ji}}$$
	for $i,j\in [n]$.
	
	Following \cite{Demazure74}, for $f\in \Bbbk[X;F]$ and $i\in [n]$, we define the {\it divided difference operator} $D_i$ on $\Bbbk[X;F]$ as
	\begin{equation*}D_{i}(f):=\frac{f-s_i(f)}{1-X_i}.
\end{equation*}
 We will need the following Leibniz rule for Demazure operators: for any $f,g\in \Bbbk(X;F)$ we have
	\begin{equation} \label{equa:Leibniz rule} D_i(fg)=D_i(f)g+s_i(f)D_i(g).
		\end{equation}
Note that $D_i^2=D_i$ \cite[Subsection 5.6]{Demazure74}.

\subsection{Affine Hecke algebras} Let $q\in \Bbbk^{\times}$ with $q\neq 1$. Following Cherednik \cite[Definition 1.1]{Cherednik91}, we give the following definition. 

\begin{definition} (a) The {\itshape Hecke algebra} $H_q$ is generated over $\Bbbk$ by $T_1,\cdots, T_n$ with the following relations
	\begin{align} 
		\label{aha relation:Tq}
		(T_i + 1)(T_i - q) &= 0,\\
		\label{aha relation:braid T}
		T_i T_jT_iT_j\cdots &= T_j T_iT_jT_i\cdots (m_{ij}\ \mbox{factors on both sides})
		\end{align}
	where $m_{ij}$ is the order of $s_is_j$, $i,j\in [n]$.
	
	(b) The {\itshape affine Hecke algebra} $\mcH_q$ is defined to be the $\Bbbk$-algebra $\Bbbk[X;F]\otimes_\Bbbk H_q$ with the following relations for all admissible indices:
	\begin{equation}
		\label{aha relation:TX}	T_i X_\alpha =X_{s_i(\alpha)}T_i+(q-1)D_i(X_\alpha).
	\end{equation}
	\end{definition}

The affine Hecke algebra $\mcH_q$ is an associative unital $\Bbbk$-algebra with $H_q$ as a subalgebra.

 By \cite[(1.6)]{Cherednik91}, we have the decomposition
 \begin{equation}\label{aha decom:Tw}\mcH_q=\oplus_{w\in \mcW}T_w\Bbbk[X;F]=\oplus_{w\in \mcW}\Bbbk[X;F]T_w. \end{equation}
 Here $$T_w=T_{i_r}\cdots T_{i_1} \quad \mbox{if} \quad w=s_{i_r}\cdots s_{i_1} \quad \mbox{for} \quad r=\ell(w)$$
 depends only on $w$ by (\ref{aha relation:braid T}).

Following \cite{Cherednik91} again or \cite[3.12(c)]{Lusztig89}, we define an algebra $\overline{\mcH}_q$ to be the associated unital $\Bbbk$-algebra generated by $\Bbbk(X;F)$ and $\mcH_q$. We have the decomposition 
\begin{equation}\label{wdecom overline}\overline{\mcH}_q=\oplus_{w\in \mcW}T_w\Bbbk(X;F)=\oplus_{w\in \mcW}\Bbbk(X;F)T_w \end{equation}	
by (\ref{aha decom:Tw}). In $\overline{\mcH}_q$, we have 
\begin{equation} \label{T_if} T_if=s_i(f)T_i+(q-1)D_i(f)\end{equation}
for any $f\in \Bbbk(X;F)$ and $i\in [n]$ (see \cite[3.12(d)]{Lusztig89}). Indeed, by \cite[3.12(a)]{Lusztig89}, we can write $f=\frac{f_1}{z}$ with $f_1\in \Bbbk[X;F]$ and $z\in \mathcal{Z}-0$, where $\mathcal{Z}$ is the center of $\mcH_q$ consisting of $\mcW$-invariant polynomials in $\Bbbk[X;F]$ (see \cite[Lemma 1.3(b)]{Cherednik91} ). Therefore, we can reduce formula (\ref{T_if}) to the case where $f\in \Bbbk[X;F]$, in which case we can use the defining relation (\ref{aha relation:TX}). 
	

\subsection{Intertwiners} Following \cite{Cherednik91}, for $i\in [n]$, we define an element $\Phi_i$ as follows:
 $$\Phi_i=T_i+(q-1)/(X_i-1)\in \overline{\mcH}_q.$$ 
 
 If $w=s_{i_1} s_{i_2}\cdots s_{i_k}\in \mcW$ is a reduced expression, then we put $\Phi_w:=\Phi_{i_1}\Phi_{i_2}\cdots \Phi_{i_k}$.

\begin{proposition} \label{prop:intertwiner for H_q} $(\mathrm{i})$ The elements $\Phi_i, i\in [n]$ satisfy the relations
	\begin{align} \label{Phi^2}	\Phi_i^2&=\tfrac{(1-qX_i)(q-X_i)}{(1-X_i)^2},\\
		\label{Phi f}	\Phi_if&=s_i(f)\Phi_i, \quad \forall f\in \Bbbk(X;F).
	\end{align}	
	
	$(\mathrm{ii})$ \cite[Proposition (b)]{Cherednik91}	The elements $\Phi_w$ have the following properties: 
\begin{align} \label{PhiwX}	\Phi_wX_{\alpha}&=X_{w(\alpha)}\Phi_w, \ w\in\mcW, \alpha\in \Sigma;\\
	\label{overline aha decom:Phiw}\overline{\mcH}_q&=\oplus_{w\in \mcW}\Phi_w\Bbbk(X;F)=\oplus_{w\in \mcW}\Bbbk(X;F)\Phi_w.
\end{align}	

$(\mathrm{iii})$ \cite[Proposition (c)]{Cherednik91} The element $\Phi_w$ does not depend on the choice of reduced decomposition of $w$:
\begin{equation}	\label{Phi braid relation}\Phi_w\Phi_{w'}=\Phi_{ww'} \ \mbox{if} \ \ell(ww')=\ell(w)+\ell(w').
 \end{equation}	
\end{proposition}
	
\begin{proof}
	It remains only to verify (i). By (\ref{T_if}) and the defining relation (\ref{aha relation:Tq}), we have
	\begin{align*}
		\Phi_i^2&=(T_i+\tfrac{q-1}{X_i-1})(T_i+\tfrac{q-1}{X_i-1})\\
		&=T_i^2+T_i\tfrac{q-1}{X_i-1}+\tfrac{q-1}{X_i-1}T_i+\tfrac{(q-1)^2}{(X_i-1)^2}\\
		&=(q-1)T_i+q+s_i(\tfrac{q-1}{X_i-1})T_i+(q-1)D_i(\tfrac{q-1}{X_i-1})+\tfrac{q-1}{X_i-1}T_i+\tfrac{(q-1)^2}{(X_i-1)^2}\\
		&=(q-1)T_i+q+\tfrac{q-1}{X_i^{-1}-1}T_i-\tfrac{(q-1)^2(X_i+1)}{(X_i-1)^2}+\tfrac{q-1}{X_i-1}T_i+\tfrac{(q-1)^2}{(X_i-1)^2}\\
		&=\tfrac{(1-qX_i)(q-X_i)}{(1-X_i)^2}.
	\end{align*}

Relation (\ref{Phi f}) follows immediately from (\ref{T_if}). 
\end{proof}

	\subsection{The modified forms of $\mcH_q$} Let $e$ be the smallest positive integer such that $1+q+\cdots +q^{e-1}=0$, setting $e:=0$ if no such integer exists.
	Let $I$ be the abelian group $\mathbb{Z}/e\mathbb{Z}$. Denote by $I^\Sigma=\{\bi=(i_\alpha) \ | \ \alpha\in \Sigma, i_\alpha\in I\}$, where $i_{\alpha+\beta}=i_\alpha+i_\beta$ for $\alpha, \beta\in \Sigma,\alpha+\beta\in \Sigma\cup\{0\}$ and $i_0=0$. Then there is an action of the Weyl group $\mcW$ on $I^\Sigma$ via
	$$w(\bi)_\alpha=i_{w^{-1}(\alpha)}, \ \bi\in I^\Sigma, w\in \mcW.$$
		
	Throughout this paper, we fix a $\mcW$-orbit $\mcC$ of $I^\Sigma$. Let $\mathcal{E}(\mcC)$ be the unital $\Bbbk$-algebra with basis $\{\epsilon(\bi) | \bi\in \mcC\}$ whose elements satisfy 
	\begin{equation} \label{deepsilon} \epsilon(\bi)\epsilon(\bj)=\delta^\bi_\bj\epsilon(\bi), \quad \sum_{\bi\in \mcC}\epsilon(\bi)=1\end{equation}
	where $\delta_{\bj}^{\bi}$ is the Kronecker symbol. Note that there is an action of $\mcW$ on $\mcE(\mcC)$ via 
	$$w(\epsilon(\bi))=\epsilon(w(\bi)), \ \bi\in \mcC, w\in \mcW.$$ 
	\begin{definition}
		\label{defn:extension of Hq} Let $\mcH_q$ be an affine Hecke algebra. The {\itshape modified form} $\overline{\mcH}_q(\mcC)$ of $\mcH_q$ associated with $\mathcal{C}$ is defined to be the $\Bbbk$-algebra $\overline{\mcH}_q(\mcC)=\oplus_{\bi\in \mcC}\overline{\mcH}_q\epsilon(\bi)$ with the relations:
		\begin{align}
			\label{Xepsi}X_\alpha\epsilon(\bi)&=\epsilon(\bi)X_\alpha, \ \bi\in \mcC, \alpha\in \Sigma\\
			\label{Teps} T_r\epsilon(\bi)&=\epsilon(s_r(\bi))T_r+\tfrac{q-1}{X_r-1}(\epsilon(s_r(\bi))-\epsilon(\bi)),\bi\in \mcC,r\in [n].		
		\end{align}
	\end{definition}
This defines an associative unital $\Bbbk$-algebra, since relation (\ref{Teps}) is equivalent to 
\begin{equation}
	\label{Phiepsi}\Phi_r\epsilon(\bi)=\epsilon(s_r(\bi))\Phi_r, \ \bi\in \mcC, r\in [n].
\end{equation}
by the decompositions (\ref{wdecom overline}) and (\ref{overline aha decom:Phiw}). Moreover, we have the decompositions 
\begin{equation} \label{extensiondecom-Tw}\overline{\mcH}_q(\mcC)=\oplus_{w\in \mcW,\bi\in \mcC}T_w\Bbbk(X;F)\epsilon(\bi)=\oplus_{w\in \mcW,\bi\in \mcC}\Bbbk(X;F)T_w\epsilon(\bi).\end{equation}			
and
\begin{equation} \label{extensiondecom:Phiw}\overline{\mcH}_q(\mcC)=\oplus_{w\in \mcW,\bi\in \mcC}\Phi_w\Bbbk(X;F)\epsilon(\bi)=\oplus_{w\in \mcW,\bi\in \mcC}\Bbbk(X;F)\Phi_w\epsilon(\bi).\end{equation}

		\subsection{The Brundan-Kleshchev auxiliary elements}
	Following \cite{Brundan-Kleshchev09}, for each $\alpha\in \Sigma$, we introduce the element 
	\begin{equation*} \label{YX} Y_\alpha:=\sum_{\bi\in \mcC}(1-q^{-i_\alpha}X_\alpha)\epsilon(\bi)\in \overline{\mcH}_q(\mcC).\end{equation*} 
	 It is easy to see that $Y_\alpha\epsilon(\bi)=(1-q^{-i_\alpha}X_\alpha)\epsilon(\bi)$ and $X_\alpha\epsilon(\bi)=q^{i_\alpha}(1-Y_\alpha)\epsilon(\bi)$. Then the elements $Y_\alpha$ satisfy $Y_\alpha\epsilon(\bi)=\epsilon(\bi)Y_\alpha$ and the relations 
$$G\colon Y_{\alpha+\beta}=Y_\alpha+Y_{\beta}-Y_{\alpha}Y_{\beta} \ \mbox{and} \ Y_0=0$$ for $\alpha, \beta\in \Sigma, \alpha+\beta\in \Sigma\cup \{0\}$. Moreover, $Y_\alpha$ is a unit in $\overline{\mcH}_q(\mcC)$ with $Y_\alpha^{-1}=\sum_{\bi\in \mcC}(1-q^{-i_\alpha}X_\alpha)^{-1}\epsilon(\bi)$. 	

Let $\Bbbk[Y;G]$ be the polynomial algebra generated by $\{Y_\alpha\ | \ \alpha\in \Sigma\}$ over $\Bbbk$ with the relations $G$. We use $\Bbbk(Y;G)$ to denote the field of fractions of $\Bbbk[Y;G]$. There is an action of the Weyl group $\mcW$ on $\Bbbk(Y;G)$ via $$w(Y_\alpha)=Y_{w(\alpha)},\quad w\in \mcW, \alpha\in \Sigma.$$	
	Note that $w(Y_\alpha):=\sum_{\bi\in \mcC}(1-q^{-i_{\alpha}}X_{w(\alpha)})\epsilon(w(\bi))$. Set $Y_r:=Y_{\alpha_r}$ for $r\in[n]$. We denote by $\parti_r$ the following Demazure operator on $\Bbbk(Y;G)$: 
	$$\parti_r(f)=\frac{s_r(f)-f}{Y_{r}}$$ 
	 which is defined in \cite{Demazure73}.

	By (\ref{Phi f}) and (\ref{Phiepsi}), the following relations hold:
	\begin{equation}\label{equa:fPhi}
	\Phi_r	f=s_r(f)\Phi_r, \quad r\in [n], f\in \Bbbk(Y;G).
	\end{equation}

For $\alpha\in \Sigma, \bi\in \mcC$, we define the element in $\overline{\mcH}_q(\mcC)$:
	\begin{equation*}\label{Qr(i)} Q_\alpha(\bi)=\begin{cases} 1-q-Y_\alpha & \text{if $i_\alpha=0$};\\
			1& \text{if $i_\alpha=1, e\neq 2$};\\
			[q-q^{-1}(1-Y_\alpha)][1-q^{-1}(1-Y_\alpha)]^{-2}& \text{if $i_\alpha=-1, e\neq 2$};\\
			[q(1-Y_\alpha)-1]^{-1}& \text{if $i_\alpha=1, e=2$};\\
			[q^{i_\alpha}(1-Y_\alpha)-q][q^{i_\alpha}(1-Y_\alpha)-1]^{-1}& \text{if $i_\alpha\neq 0, \pm 1$.}
		\end{cases}
	\end{equation*}
	We set $Q_r(\bi):=Q_{\alpha_r}(\bi)$ and $i_{r}:=i_{\alpha_r}$ for simple root $\alpha_r$ and $\bi\in \mcC$. For $w\in \mcW, \alpha\in \Sigma $ and $\bi\in \mcC$, it is easy to check that $$w(Q_\alpha(\bi))=Q_{w(\alpha)}(w(\bi))$$ and 
		\begin{equation}\label{equa:QPhi}
		\Phi_rQ_\alpha(\bi)=Q_{s_r(\alpha)}(s_r(\bi))\Phi_r, \ r\in [n]
	\end{equation}
 by (\ref{equa:fPhi}).
 
For $\alpha\in \Sigma$ and $\bi\in \mcC$, note that $Y_{-\alpha}=Y_\alpha(Y_\alpha-1)^{-1}$. A direct computation shows that the element $Q_\alpha(\bi)$ satisfies the following formula:
			\begin{equation}\label{QrsrQrsr} Q_{s_\alpha(\alpha)}(\bi)Q_\alpha(\bi)=\begin{cases}\tfrac{[(1-Y_\alpha)-q][1-q(1-Y_\alpha)]}{1-Y_\alpha}& \text{if $i_\alpha=0$};\\
				\tfrac{(1-Y_\alpha)[q(1-Y_\alpha)-q^{-1}]}{(1-Y_\alpha-q^{-1})^2}&\text{if $i_\alpha=1, e\neq 2$};\\
				\tfrac{q-q^{-1}(1-Y_\alpha)}{[1-q^{-1}(1-Y_\alpha)]^2}&\text{if $i_\alpha=-1, e\neq2$};\\
				\tfrac{q(Y_\alpha-1)}{[1-q(1-Y_\alpha)]^2}&\text{if $i_\alpha=1, e=2$};\\
				\tfrac{[1-q^{i_\alpha+1}(1-Y_\alpha)][q-q^{i_\alpha}(1-Y_\alpha)]}{[1-q^{i_\alpha}(1-Y_\alpha)]^2}& \text{if $i_\alpha\neq 0, \pm 1$}.
			\end{cases}
		\end{equation}

	For each $\alpha\in \Sigma$ and $\bi\in \mcC$, inside $\overline{\mcH}_q(\mcC)$, we denote by 
	\begin{equation}\label{P_r} P_\alpha(\bi)=\begin{cases}Y_\alpha^{-1}-Y_{\alpha}^{-2}& \text{if $i_\alpha=0$};\\
			Y_{\alpha}(Y_\alpha-1) &\text{if $i_\alpha=1, e\neq 2$};\\
			Y_\alpha & \text{if $i_\alpha=-1, e\neq 2$};\\
			Y_\alpha^2 (Y_{\alpha}-1)^{-1}& \text{if $i_\alpha=1, e=2$};\\
			1& \text{if $i_\alpha\neq0,\pm1$},
		\end{cases}
	\end{equation}	
	We set $P_r(\bi):=P_{\alpha_r}(\bi)$ for simple root $\alpha_r$ and $\bi\in \mcC$. Notice that 
	$$w(P_\alpha(\bi))=P_{w(\alpha)}(w(\bi)).$$

	For each $r\in [n]$, we define $$\Theta_r=\Phi_r\sum_{\bi\in \mcC}Q_r^{-1}(\bi)\epsilon(\bi)\in \overline{\mcH}_q(\mcC).$$
	If $w=s_{r_1} s_{r_2}\cdots s_{r_k}\in \mcW$ is a reduced expression, then we put $\Theta_w:=\Theta_{r_1}\Theta_{r_2}\cdots \Theta_{r_k}$. The following result is a modification of Proposition \ref{prop:intertwiner for H_q}.
	\begin{proposition}
		\label{lem:Theta} $(\mathrm{i})$ The elements $\Theta_r$ satisfy the following relations:
		\begin{align} \label{Thetaepsi}\Theta_r\epsilon(\bi)&=\epsilon(s_r(\bi))\Theta_r, \quad \forall r\in [n], \bi\in \mcC; \\
			\label{fTheta} \Theta_rf&=s_r(f)\Theta_r, \quad \forall r\in [n], f\in \Bbbk(Y;G).
		\end{align}			
		
		$(\mathrm{ii})$	For $r\in[n], \bi\in \mcC$, we have:		
		\begin{equation}\label{Theta^2} \Theta_r^2\epsilon(\bi)=P_r(\bi)\epsilon(\bi).		\end{equation}	
		
	$(\mathrm{iii})$	The elements $\Theta_w$ have the following properties:
				\begin{align} \label{ThetawY}	\Theta_wY_{\alpha}&=Y_{w(\alpha)}\Theta_w, \ \forall w\in\mcW, \alpha\in \Sigma,\\
			\label{ThetaoverlineH(E)}\overline{\mcH}_q(\mcC)&=\oplus_{w\in \mcW,\bi\in \mcC}\Theta_w\Bbbk(Y;G)\epsilon(\bi)=\oplus_{w\in \mcW,\bi\in \mcC}\Bbbk(Y;G)\Theta_w\epsilon(\bi).
		\end{align}

		$(\mathrm{iv})$ The element $\Theta_w$ does not depend on the choice of reduced decomposition of $w$:
		\begin{equation} \label{Theta braid relation} \Theta_w\Theta_{w'}=\Theta_{ww'} \quad \mbox{if} \ \ell(ww')=\ell(w)+\ell(w').
			\end{equation}
		
	\end{proposition}
	
	\begin{proof}
		(i) Relation (\ref{Thetaepsi}) follows from the commutativity of the elements $Y_\alpha$ and $\epsilon(\bi)$ and the relation (\ref{Phiepsi}). Relation (\ref{fTheta}) follows from (\ref{equa:fPhi}) and (\ref{Thetaepsi}) .
		
		(ii) By (\ref{Phiepsi}) and (\ref{equa:QPhi}), we have
		\begin{align*}
			\Theta_r^2\epsilon(\bi)&=(\Phi_r\sum_{\bj\in \mcC}Q^{-1}_r(\bj)\epsilon(\bj))\Phi_rQ_r^{-1}(\bi)\epsilon(\bi)\\
			&=\Phi^2_r\sum_{\bj\in \mcC}Q^{-1}_{s_r(\alpha_r)}(s_r(\bj))\epsilon(s_r(\bj))Q^{-1}_r(\bi)\epsilon(\bi)\\
			&=\Phi^2_rQ_{s_r(\alpha_r)}^{-1}(\bi)Q_{r}(\bi)^{-1}\epsilon(\bi).
		\end{align*}
Combining this identity with (\ref{Phi^2}) and (\ref{QrsrQrsr}) gives (\ref{Theta^2}).
	
		(iii) Relation (\ref{ThetawY}) follows from (\ref{fTheta}) and the definition of $\Theta_w$. Applying (\ref{equa:fPhi}) repeatedly we can write
			\begin{equation}\label{equa:Theta-Phi}\Theta_w\epsilon(\bi)=\Phi_w g\epsilon(\bi)\end{equation}
		for some $g\in \Bbbk(Y;G)$. Moreover, $\Bbbk(X;F)\epsilon(\bi)=\Bbbk(Y;G)\epsilon(\bi)$ as $\Bbbk$-vector spaces for each $\bi\in \mcC$. 
		Therefore the decomposition (\ref{ThetaoverlineH(E)}) follows from (\ref{extensiondecom:Phiw}).
		
		(iv) Assume $w=s_{r_1} s_{r_2}\cdots s_{r_k}\in \mcW$ is a reduced expression. For each $\bi\in \mcC$, relations (\ref{Thetaepsi}) and (\ref{equa:QPhi}) give 
		\begin{align*}\Theta_w\epsilon(\bi)&=\Theta_{r_1}\cdots \Theta_{r_{k-1}}\Theta_{r_k}\epsilon(\bi)\\
			&			=\Phi_wQ^{-1}_{s_{r_k}\cdots s_{r_2}(\alpha_{r_1})}(\bi)\cdots Q^{-1}_{s_{r_k}(\alpha_{r_{k-1}})}(\bi)Q^{-1}_{r_k}(\bi)\epsilon(\bi).
		\end{align*}
		Thus, applying (\ref{Phi braid relation}), we know that the elements $\Theta_r$ satisfy the Coxeter relations of $\mcW$:
	$$	\Theta_{r}\Theta_{t}\Theta_{r}\cdots=\Theta_{t}\Theta_{r}\Theta_{t}\cdots (m_{rt}\ \mbox{factors on both sides}),\ r,t\in [n].$$
	Therefore, the element $\Theta_w$ does not depend on the choice of reduced decomposition of $w$ by \cite[Theorem 3.3.1(ii)]{Bjorner-Brenti}. In particular, we have formula (\ref{Theta braid relation}).	\end{proof}

	\section{KLR-like generators of $\overline{\mcH}_q(\mcC)$}
	
	We retain the notation of Section~2. 
	For each $r\in [n]$, we define an element	
	 \begin{equation} \label{Psitheta}\Psi_r=\sum_{\bi\in \mcC}(\Theta_r-\delta_{i_r}^{0}Y_r^{-1})\epsilon(\bi)\in \overline{\mcH}_q(\mcC).
	\end{equation}
The aim of this section is to prove the following result, which gives a KLR-like presentation of $\overline{\mcH}_q(\mcC)$.

	\begin{theorem} \label{thm:KLR basis} The algebra $\overline{\mcH}_q(\mcC)$ is generated by 
	$$\{\Psi_r\ | \ r\in [n]\}\cup \{\epsilon(\bi)\ | \ \bi\in \mcC\}\cup \{f \ | \ f\in \Bbbk(Y;G)\}$$
		 and the following is a set of complete relations for all $\alpha\in \Sigma, \bi\in \mcC$ and $r,t\in [n]$.
		\begin{enumerate}
			
			\item $\epsilon(\bi)\epsilon(\bj)=\delta_{\bi}^{\bj}\epsilon(\bi), \ \sum_{\bi\in \mcC}\epsilon(\bi)=1$.
			
			\item $Y_\alpha\epsilon(\bi)=\epsilon(\bi)Y_\alpha$.			
			
		\item	$\Psi_r\epsilon(\bi)=\epsilon(s_r(\bi))\Psi_r$.			
			
	\item		 $\Psi_rY_\alpha\epsilon(\bi)=[Y_{s_r(\alpha)}\Psi_r+\delta_{i_r}^{0}\parti_r(Y_\alpha)]\epsilon(\bi)$.
			
	\item $\Psi_r^2\epsilon(\bi)=	\begin{cases}-\Psi_r\epsilon(\bi)&\text{if $i_r=0$};\\
			Y_{-\alpha_r}\epsilon(\bi) & \text{if $i_r=1, e\neq 2$};\\
		Y_{r}\epsilon(\bi) & \text{if $i_r=-1, e\neq 2$};\\
		Y_{r}Y_{-\alpha_r}\epsilon(\bi) & \text{if $i_r=1, e=2$};\\
		\epsilon(\bi)&\text{if $i_r\neq0,\pm 1$}.
	\end{cases}$
	
	\item If $r\nslash t$ then 
	$\Psi_r\Psi_t=\Psi_t\Psi_r.$
	
	\item If $r-\-t$ then {\tiny $$(\Psi_r\Psi_t\Psi_r-\Psi_t\Psi_r\Psi_t)\epsilon(\bi)=\begin{cases}(1-Y_{-\alpha_r})\epsilon(\bi) &\text{if $i_r=-i_t=1,e\neq 2$};\\
			(Y_{-\alpha_t}-1)\epsilon(\bi) &\text{if $i_r=-i_t=-1,e\neq 2$};\\
			(Y_{-\alpha_t}-Y_{-\alpha_r})\epsilon(\bi) &\text{if $i_r=-i_t=1,e=2$};\\
			0 & \text{else}.
		\end{cases}$$}
	
	\item If $\xymatrix@C10pt{r\ar@{=}[r]|{\rangle}\ar@{=}[r]&t}$ then {\tiny
		\begin{align*}&(\Psi_r\Psi_t\Psi_r\Psi_t-\Psi_t\Psi_r\Psi_t\Psi_r)\epsilon(\bi)\\
			&=\begin{cases}(1-Y_{-\alpha_r-\alpha_t})\Psi_r\epsilon(\bi) &\text{if $i_r=-2,i_t=1,e\neq 2$};\\
				(Y_{-\alpha_t}-1)\Psi_r\epsilon(\bi) &\text{if $i_r=2,i_t=-1,e\neq 2$};\\
				(1-Y_{-\alpha_r-\alpha_t})(Y_r\Psi_r+1)\epsilon(\bi) &\text{if $i_r=0,i_t=1,e=2$};\\
				(Y_{-\alpha_t}+Y_{-\alpha_r-2\alpha_t}-2)\Psi_t\epsilon(\bi) &\text{if $i_r=-i_t=1,e\neq 2$};\\
				(2-Y_t-Y_{-\alpha_r})\Psi_t\epsilon(\bi) &\text{if $i_r=-i_t=-1,e\neq 2$};\\
				(Y_{-\alpha_t}+Y_{-\alpha_r-2\alpha_t}-Y_{-\alpha_r}-Y_t)\Psi_t\epsilon(\bi) &\text{if $i_r=-i_t=1,e=2$};\\
				0 & \text{else}.
			\end{cases}
	\end{align*}}

\item If $\xymatrix@C10pt{r\ar@3{-}[r]|{\rangle}\ar@3{-}[r]&t}$ then {\tiny
	\begin{align*}\label{Psi6}
		&(\Psi_r\Psi_t\Psi_r\Psi_t\Psi_r\Psi_t-\Psi_t\Psi_r\Psi_t\Psi_r\Psi_t\Psi_r)\epsilon(\bi)\\
		&=\begin{cases} (1-Y_{-\alpha_r-2\alpha_t})\Psi_{r}\Psi_t\Psi_r\epsilon(\bi) &\text{if $i_r=-3, i_t=1, e\neq 2,3$};\\
			(Y_{-\alpha_t}-1)\Psi_{r}\Psi_t\Psi_r\epsilon(\bi) &\text{if $i_r=3, i_t=-1, e\neq 2,3$};\\	
			(Y_{-\alpha_t}+Y_{-\alpha_r-2\alpha_t}+Y_{-2\alpha_r-3\alpha_t}-3)\Psi_{t}\Psi_r\Psi_t\epsilon(\bi) &\text{if $i_r=-i_t=1, e\neq 2$};\\
			(3-Y_t-Y_{-\alpha_r}-Y_{\alpha_r+2\alpha_t})\Psi_{t}\Psi_r\Psi_t\epsilon(\bi) &\text{if $i_r=-i_t=-1, e\neq 2$};\\
			(Y_{\alpha_r+2\alpha_t}-2)(Y_{-\alpha_r-\alpha_t}Y_{-2\alpha_r-3\alpha_t}-Y_{-\alpha_r})	\Psi_{r}\epsilon(\bi) &\text{if $i_r=-2, i_t=1, e\neq 2$};\\
			(2-Y_{\alpha_r+2\alpha_t})(Y_{-\alpha_t}Y_{-\alpha_r-3\alpha_t}-Y_r)\Psi_{r}\epsilon(\bi) &\text{if $i_r=2, i_t=-1, e\neq 2$};\\	
			(2-Y_{t})[Y_{-\alpha_r-2\alpha_t}(Y_{-\alpha_r-3\alpha_t}+Y_r)-Y_rY_{-2\alpha_r-3\alpha_t}](Y_t\Psi_{t}+1)\epsilon(\bi)&\text{if $i_r=1,i_{t}=0, e=2$},\\
			(Y_{-\alpha_r-\alpha_t}-Y_{-\alpha_r-3\alpha_t})	\Psi_{t}\epsilon(\bi) &\text{if $i_r=1,i_{t}=2, e=4$},\\
			(Y_{-\alpha_r}-Y_{-\alpha_r-2\alpha_t})	\Psi_{t}\epsilon(\bi) &\text{if $i_r=-1,i_{t}=2, e=4$},\\
			(Y_{-\alpha_r-3\alpha_t}-1)\Psi_t\epsilon(\bi)& \text{if $i_r=1,3i_{t}=-2, e\neq 2, 4$},\\
			(1-Y_{-\alpha_r})\Psi_t\epsilon(\bi)& \text{if $i_r=-1, 3i_{t}=2, e\neq 2, 4$},\\
			(Y_{-\alpha_r-2\alpha_t}-1)\Psi_t\epsilon(\bi) &\text{if $i_r=3,i_{t}=-2, e\neq 2,3,4$},\\
			(1-Y_{-\alpha_r-\alpha_t})\Psi_t\epsilon(\bi) &\text{if $i_r=-3,i_{t}=2, e\neq 2,3,4$},\\
			(Y_{\alpha_r+2\alpha_t}-2)[Y_{-2\alpha_r-3\alpha_t}(Y_{-\alpha_r-\alpha_t}+Y_{-\alpha_t})-Y_{-\alpha_r}Y_{-\alpha_t}](Y_r\Psi_{r}+1)\epsilon(\bi) &\text{if $i_r=0, i_{t}=1, e=2$}, \\
			[2Y_{-\alpha_t}Y_{-2\alpha_r-3\alpha_t}-Y_{-\alpha_r-\alpha_t}Y_{-\alpha_r-3\alpha_t}+(Y_{-\alpha_t}-Y_{-\alpha_r-2\alpha_t})\Psi_r\Psi_t\Psi_r	\\
		+(Y_{-\alpha_t}+Y_{-\alpha_r-2\alpha_t}+Y_{-2\alpha_r-3\alpha_t}-Y_{-\alpha_r}-Y_t-Y_{\alpha_r+2\alpha_t})\Psi_{t}\Psi_r\Psi_t	]\epsilon(\bi) &\text{if $i_r=i_{t}=1, e=2$},\\
			(1-Y_{-\alpha_r-2\alpha_t})(\Psi_{r}\Psi_t\Psi_r+\Psi_r\Psi_t)\epsilon(\bi) &\text{if $i_r=0, i_t=1, e=3$};\\
			(Y_{-\alpha_t}-1)\Psi_{r}\Psi_t\Psi_r\epsilon(\bi) &\text{if $i_r=0, i_t=-1, e=3$};\\	
			0&\text{else.}
		\end{cases}	
\end{align*}}
		\end{enumerate}
		\end{theorem}
	
For readability, we divide the verification of the relations of Theorem~\ref{thm:KLR basis} into several lemmas. Statements (1) and (2) follow directly from the definitions of $\epsilon(\bi)$ and $Y_\alpha$. 
	
	\begin{lemma}
		\label{lem:Psiepsilon} For all $r\in [n]$ and $\bi\in \mcC$, we have 
		$$\Psi_r\epsilon(\bi)=\epsilon(s_r(\bi))\Psi_r.$$
	\end{lemma}
	\begin{proof} By (\ref{Thetaepsi}), we get that
		\begin{equation*}\Psi_r\epsilon(\bi)=(\Theta_r-\delta_{i_r}^0Y_r^{-1})\epsilon(\bi)=\epsilon(s_r(\bi))(\Theta_r-\delta_{i_r}^0Y_r^{-1})=\epsilon(s_r(\bi))\Psi_r\end{equation*}
		where the third identity holds since $s_r(\bi)=\bi$ whenever $i_r=0$.
		\end{proof}
	
	We first record a more general form of relation (4) in Theorem~\ref{thm:KLR basis}.
	\begin{lemma}
		\label{lem:PsiY} For all $r\in [n]$ and $\bi\in \mcC$, we have 
		$$\Psi_r f\epsilon(\bi)=(s_r(f)\Psi_r+\delta^0_{i_r}\parti_r(f))\epsilon(\bi), \ \forall f\in \Bbbk(Y;G).$$
	\end{lemma}
	\begin{proof}
		By (\ref{fTheta}),
		\begin{align*}\Psi_rf\epsilon(\bi)&=(\Theta_r-\delta_{i_r}^0Y_r^{-1})f\epsilon(\bi)\\
			&=(s_r(f)\Theta_r-\delta_{i_r}^0f Y_r^{-1})\epsilon(\bi) \\
			&=[s_r(f)(\Psi_r+\delta_{i_r}^0Y_r^{-1})-\delta_{i_r}^0f Y_r^{-1}]\epsilon(\bi)\\
			&=(s_r(f)\Psi_r+\delta_{i_r}^0\parti_r(f))\epsilon(\bi).
		\end{align*} 
	\end{proof}

\begin{lemma}
	\label{lem:Psi^2} For all $r\in [n]$ and $\bi\in \mcC$, we have
	$$\Psi_r^2\epsilon(\bi)=	\begin{cases}-\Psi_r\epsilon(\bi)&\text{if $i_r=0$};\\
		Y_{-\alpha_r}\epsilon(\bi) & \text{if $i_r=1, e\neq 2$};\\
		Y_{r}\epsilon(\bi) & \text{if $i_r=-1, e\neq 2$};\\
		Y_{r}Y_{-\alpha_r}\epsilon(\bi) & \text{if $i_r=1, e=2$};\\
		\epsilon(\bi)&\text{if $i_r\neq0,\pm 1$}.
		\end{cases}$$
\end{lemma}

\begin{proof}
	 Using Lemma~\ref{lem:Psiepsilon} and (\ref{Thetaepsi}), and noting that $s_r(\bi)_r=-i_r$, we obtain
	\begin{align*}
		\Psi_r^2\epsilon(\bi)&=\Psi_r\epsilon(s_r(\bi))\Psi_r\epsilon(\bi)\\
		&=(\Theta_r-\delta_{i_r}^0Y_r^{-1})(\Theta_r-\delta_{i_r}^0Y_r^{-1})\epsilon(\bi)\\
		&=(\Theta_r^2-\delta_{i_r}^0Y_{s_r(\alpha_r)}^{-1}\Theta_r-\delta_{i_r}^0Y_r^{-1}\Theta_r+\delta_{i_r}^0Y^{-2}_r)\epsilon(\bi)\\
		&=[\Theta_r^2-\delta_{i_r}^0Y_{r}^{-1}(Y_r-1)\Theta_r-\delta^0_{i_r}Y_r^{-1}\Theta_r+\delta_{i_r}^0Y^{-2}_r]\epsilon(\bi)\\
		&=[\Theta_r^2-\delta_{i_r}^0\Theta_r+\delta_{i_r}^0Y^{-2}_r]\epsilon(\bi)\\
		&=\begin{cases}(Y_r^{-1}-Y_{r}^{-2}-\Theta_r+Y_r^{-2})\epsilon(\bi)&\text{if $i_r=0$};\\
			P_r(\bi)&\text{if $i_r\neq 0$}
		\end{cases}\\
	&=\begin{cases}-\Psi_r\epsilon(\bi)&\text{if $i_r=0$};\\
		P_r(\bi)&\text{if $i_r\neq 0$}.
	\end{cases}
	\end{align*}
Then the claim follows from the fact that $Y_{-\alpha}=Y_\alpha(Y_\alpha-1)^{-1}$ for all $\alpha\in \Sigma$.
\end{proof}

\begin{lemma}
	\label{lem:r nslash t} Let $r,t\in [n]$. If $r\nslash t$, then 
	$$\Psi_r\Psi_t=\Psi_t\Psi_r.$$
\end{lemma}
\begin{proof}
	 Since $r\nslash t$, by Lemma \ref{lem:Psiepsilon}, (\ref{Thetaepsi}) and (\ref{Theta braid relation}), for any $\bi\in \mcC$ we have 
	\begin{align*}\Psi_r\Psi_t\epsilon(\bi)&=\Psi_r\epsilon(s_t(\bi))\Psi_t\epsilon(\bi)\\
		&=(\Theta_r-\delta_{i_r}^0Y_r^{-1})(\Theta_t-\delta_{i_t}^0Y_t^{-1})\epsilon(\bi)\\
		&=(\Theta_{s_rs_t}-\delta_{i_r}^0Y_r^{-1}\Theta_t-\delta_{i_t}^0\Theta_rY_t^{-1}+\delta_{i_r}^0\delta_{i_t}^0Y_r^{-1}Y_t^{-1})\epsilon(\bi)\\
		&=(\Theta_{s_ts_r}-\delta_{i_r}^0\Theta_tY_r^{-1}-\delta_{i_t}^0Y_t^{-1}\Theta_r+\delta_{i_r}^0\delta_{i_t}^0Y_r^{-1}Y_t^{-1})\epsilon(\bi)\\
		&=(\Theta_t-\delta_{i_t}^0Y_t^{-1})(\Theta_r-\delta_{i_r}^0Y_r^{-1})\epsilon(\bi)	\\
		&=\Psi_t\Psi_r\epsilon(\bi).
	\end{align*}
	Thus we have $\Psi_r\Psi_t=\Psi_t\Psi_r$.\end{proof}

\begin{lemma}
	\label{lem: r--t} Let $r,t\in [n]$. If $r-\-t$ then Theorem \ref{thm:KLR basis} (7) holds.
\end{lemma}

\begin{proof}
	 In this case, by Lemma \ref{lem:Psiepsilon} and (\ref{Thetaepsi})-(\ref{Theta^2}), we have:
	\begin{align*}\Psi_r\Psi_t\Psi_r\epsilon(\bi)&=\Psi_r\epsilon(s_ts_r(\bi))\Psi_t\epsilon(s_r(\bi))\Psi_r\epsilon(\bi)\\
		&=(\Theta_r-\tfrac{\delta_{s_ts_r(\bi)_r}^0}{Y_r})(\Theta_t-\tfrac{\delta_{s_r(\bi)_t}^0}{Y_t})(\Theta_r-\tfrac{\delta_{i_r}^0}{Y_r})\epsilon(\bi)\\
		&=(\Theta_r-\tfrac{\delta_{i_t}^0}{Y_r})(\Theta_t-\tfrac{\delta_{i_r+i_t}^0}{Y_t})(\Theta_r-\tfrac{\delta_{i_r}^0}{Y_r})\epsilon(\bi)\\
		&=(\Theta_{s_rs_t}-\tfrac{\delta_{i_r+i_t}^{0}}{Y_{s_r(\alpha_t)}}\Theta_r-\tfrac{\delta_{i_t}^0}{Y_r}\Theta_t+\tfrac{\delta_{i_r}^0\delta_{i_t}^0}{Y_rY_t})(\Theta_r-\tfrac{\delta_{i_r}^0}{Y_r})\epsilon(\bi)\\
		&=[\Theta_{s_rs_ts_r}-\tfrac{\delta_{i_r}^0}{Y_{s_rs_t(\alpha_r)}}\Theta_{s_rs_t}-\tfrac{\delta_{i_{t}}^{0}}{Y_r}\Theta_{s_ts_r}+(\tfrac{\delta_{i_r}^0\delta_{i_{t}}^{0}}{Y_{s_r(\alpha_t)}Y_{s_r(\alpha_r)}}+\tfrac{\delta_{i_r}^0\delta_{i_{t}}^{0}}{Y_{r}Y_t})\Theta_r\\
		&+\tfrac{\delta_{i_r}^0\delta_{i_{t}}^{0}}{Y_{\alpha_r+\alpha_t}Y_r}\Theta_t-\tfrac{\delta_{i_{r}+i_t}^0P_r(\bi)}{Y_{s_r(\alpha_t)}}-\tfrac{\delta_{i_r}^0\delta_{i_{t}}^{0}}{Y_r^2Y_t}]\epsilon(\bi)\\
		&=[\Theta_{s_rs_ts_r}-\tfrac{\delta_{i_r}^0}{Y_{t}}\Theta_{s_rs_t}-\tfrac{\delta_{i_{t}}^{0}}{Y_r}\Theta_{s_ts_r}+\tfrac{\delta_{i_r}^0\delta_{i_{t}}^{0}}{Y_{\alpha_r+\alpha_t}Y_t}\Theta_r+\tfrac{\delta_{i_r}^0\delta_{i_{t}}^{0}}{Y_{\alpha_r+\alpha_t}Y_r}\Theta_t\\
		&-\tfrac{\delta_{i_{r}+i_t}^0P_r(\bi)}{Y_{\alpha_r+\alpha_t}}-\tfrac{\delta_{i_r}^0\delta_{i_{t}}^{0}}{Y_r^2Y_t}]\epsilon(\bi)
	\end{align*}
	Similarly, we can show that
	\begin{align*}&\Psi_t\Psi_r\Psi_t\epsilon(\bi)\\
		&=[\Theta_{s_ts_rs_t}-\tfrac{\delta_{i_t}^0}{Y_{r}}\Theta_{s_ts_r}-\tfrac{\delta_{i_{r}}^{0}}{Y_t}\Theta_{s_rs_t}+\tfrac{\delta_{i_r}^0\delta_{i_{t}}^{0}}{Y_{\alpha_r+\alpha_t}Y_t}\Theta_r+\tfrac{\delta_{i_r}^0\delta_{i_{t}}^{0}}{Y_{\alpha_r+\alpha_t}Y_r}\Theta_t-\tfrac{\delta_{i_{r}+i_t}^0P_t(\bi)}{Y_{\alpha_r+\alpha_t}}-\tfrac{\delta_{i_r}^0\delta_{i_{t}}^{0}}{Y_rY_t^2}]\epsilon(\bi)
	\end{align*}
	By (\ref{Theta braid relation}), we arrive at 
	\begin{align*}
		(\Psi_r\Psi_t\Psi_r-\Psi_t\Psi_r\Psi_t)\epsilon(\bi)&=(\tfrac{\delta_{i_{r}+i_t}^0(P_t(\bi)-P_r(\bi))}{Y_{\alpha_r+\alpha_t}}+\tfrac{\delta_{i_r}^0\delta^0_{i_t}(Y_r-Y_t)}{Y_r^2Y_t^2})\epsilon(\bi).
	\end{align*}
	There are two cases:
	
	$Case~1\colon i_r=i_t=0$. Using (\ref{P_r}), we have
$$
		(\Psi_r\Psi_t\Psi_r-\Psi_t\Psi_r\Psi_t)\epsilon(\bi)=(\tfrac{(Y_t^{-1}-Y_t^{-2})-(Y_r^{-1}-Y_r^{-2})}{Y_{r}+Y_t-Y_rY_t}+\tfrac{Y_r-Y_t}{Y_r^{2}Y_t^{2}})\epsilon(\bi)=0.
	$$
	$Case~2\colon i_r=-i_t\neq 0$. Using (\ref{P_r}) again, we have
	\begin{align*}(\Psi_r\Psi_t\Psi_r-\Psi_t\Psi_r\Psi_t)\epsilon(\bi)&=\tfrac{P_t(\bi)-P_r(\bi)}{Y_{\alpha_r+\alpha_t}}\epsilon(\bi)\\
				&=\begin{cases}\tfrac{Y_t-Y_r(Y_r-1)^{-1}}{Y_r+Y_t-Y_rY_t}\epsilon(\bi) &\text{if $i_r=-i_t=1,e\neq 2$};\\
			\tfrac{Y_t(Y_t-1)^{-1}-Y_r}{Y_r+Y_t-Y_rY_t}\epsilon(\bi) &\text{if $i_r=-i_t=-1,e\neq 2$};\\
			\tfrac{Y_t^2(Y_t-1)^{-1}-Y_r^2(Y_r-1)^2}{Y_{r}+Y_t-Y_rY_t}\epsilon(\bi) &\text{if $i_r=-i_t=1,e=2$};\\
			0 & \text{else}
		\end{cases}\\
		&=\begin{cases}\tfrac{1}{1-Y_r}\epsilon(\bi) &\text{if $i_r=-i_t=1,e\neq 2$};\\
			\tfrac{1}{Y_t-1}\epsilon(\bi) &\text{if $i_r=-i_t=-1,e\neq 2$};\\
			\tfrac{Y_r-Y_t}{(Y_{r}-1)(Y_{t}-1)}\epsilon(\bi) &\text{if $i_r=-i_t=1,e=2$};\\
			0 & \text{else}
		\end{cases}\\
		&=\begin{cases}(1-Y_{-\alpha_r})\epsilon(\bi) &\text{if $i_r=-i_t=1,e\neq 2$};\\
			(Y_{-\alpha_t}-1)\epsilon(\bi) &\text{if $i_r=-i_t=-1,e\neq 2$};\\
			(Y_{-\alpha_t}-Y_{-\alpha_r})\epsilon(\bi) &\text{if $i_r=-i_t=1,e=2$};\\
			0 & \text{else}.
		\end{cases}
	\end{align*}
This proves relation (7) of Theorem~\ref{thm:KLR basis}.
	\end{proof}

\begin{lemma}
	\label{lem:r2->t} Let $r,t\in [n]$. If $\xymatrix@C10pt{r\ar@{=}[r]|{\rangle}\ar@{=}[r]&t}$, then the statement (8) of Theorem \ref{thm:KLR basis} holds.
\end{lemma}
\begin{proof}
	 Since $\xymatrix@C10pt{r\ar@{=}[r]|{\rangle}\ar@{=}[r]&t}$, by Lemma \ref{lem:Psiepsilon} and (\ref{Thetaepsi})-(\ref{Theta^2}), we deduce that
	\begin{align*}
		&\Psi_r\Psi_t\Psi_r\Psi_t\epsilon(\bi)\\
		&=\Psi_r\epsilon(s_ts_rs_t(\bi))\Psi_t\epsilon(s_rs_t(\bi))\Psi_r\epsilon(s_t(\bi))\Psi_t\epsilon(\bi)\\
		&=(\Theta_r-\tfrac{\delta_{s_ts_rs_t(\bi)_r}^0}{Y_r})(\Theta_t-\tfrac{\delta_{s_rs_t(\bi)_t}^0}{Y_t})(\Theta_r-\tfrac{\delta_{s_t(\bi)_r}^0}{Y_r})(\Theta_t-\tfrac{\delta_{i_{t}}^0}{Y_t})\epsilon(\bi)\\
		&=(\Theta_r-\tfrac{\delta_{i_r}^0}{Y_r})(\Theta_t-\tfrac{\delta_{i_r+i_t}^0}{Y_t})(\Theta_r-\tfrac{\delta_{i_r+2i_t}^0}{Y_r})(\Theta_t-\tfrac{\delta_{i_{t}}^0}{Y_t})\epsilon(\bi)\\
		&=(\Theta_{s_rs_t}-\tfrac{\delta_{i_r+i_t}^0}{Y_{s_r(\alpha_t)}}\Theta_r-\tfrac{\delta_{i_r}^0}{Y_r}\Theta_t+\tfrac{\delta_{i_r}^0\delta_{i_t}^0}{Y_rY_t})(\Theta_{s_rs_t}-\tfrac{\delta_{i_t}^0}{Y_{s_r(\alpha_t)}}\Theta_r-\tfrac{\delta_{i_r+2i_t}^0}{Y_r}\Theta_t+\tfrac{\delta_{i_r}^0\delta_{i_t}^0}{Y_rY_t})\epsilon(\bi)\\
		&=[\Theta_{s_rs_ts_rs_t}-\tfrac{\delta_{i_t}^0}{Y_{s_rs_ts_r(\alpha_t)}}\Theta_{s_rs_ts_r}-\tfrac{\delta_{i_r}^0}{Y_r}\Theta_{s_ts_rs_t}+(\tfrac{\delta_{i_r}^0\delta_{i_t}^0}{Y_{s_rs_t(\alpha_r)}Y_{s_rs_t(\alpha_t)}}+\tfrac{\delta_{i_r+i_t}^0\delta_{i_r+2i_t}^0}{Y_{s_r(\alpha_t)}Y_{s_r(\alpha_r)}}\\
		&+\tfrac{\delta_{i_r}^0\delta_{i_t}^0}{Y_rY_t})\Theta_{s_rs_t}+\tfrac{\delta_{i_r}^0\delta_{i_t}^0}{Y_rY_{s_ts_r(\alpha_t)}}\Theta_{s_ts_r}-(\tfrac{\delta_{i_r+2i_t}^0P_{s_r(\alpha_t)(s_r(\bi))}}{Y_{s_rs_t(\alpha_r)}}+\tfrac{\delta_{i_r}^0\delta_{i_t}^0}{Y_{s_r(\alpha_r)}Y_{s_r(\alpha_t)}^2}+\tfrac{\delta_{i_r}^0\delta_{i_t}^0}{Y_rY_tY_{s_r(\alpha_t)}})\Theta_r\\
		&-(\tfrac{\delta_{i_r+i_t}^0P_{r}(s_t(\bi))}{Y_{s_r(\alpha_t)}}+\tfrac{\delta_{i_r}^0\delta_{i_t}^0}{Y_rY_{s_t(\alpha_r)}Y_{s_t(\alpha_t)}}+\tfrac{\delta_{i_r}^0\delta_{i_t}^0}{Y_r^2Y_t})\Theta_t+\tfrac{\delta_{i_r+i_t}^0\delta_{i_t}^0P_r(\bi)}{Y_tY_{s_r(\alpha_t)}}+\tfrac{\delta_{i_r}^0\delta_{i_r+2i_t}^0P_t(\bi)}{Y_rY_{s_t(\alpha_r)}}+\tfrac{\delta_{i_r}^0\delta_{i_t}^0}{Y_r^2Y_t^2}]\epsilon(\bi)\\
		&=[\Theta_{s_rs_ts_rs_t}-\tfrac{\delta_{i_t}^0}{Y_t}\Theta_{s_rs_ts_r}-\tfrac{\delta_{i_{r}}^0}{Y_r}\Theta_{s_ts_rs_t}+\tfrac{\delta_{i_{r}}^0\delta_{i_{t}}^0}{Y_tY_{\alpha_r+2\alpha_t}}\Theta_{s_rs_t}+\tfrac{\delta_{i_{r}}^0\delta_{i_t}^0}{Y_rY_{\alpha_r+\alpha_t}}\Theta_{s_ts_r}\\
		&-(\tfrac{\delta_{i_{r}}^0\delta_{i_t}^0}{Y_tY_{\alpha_r+\alpha_t}^2}+\tfrac{\delta_{i_{r}+2i_t}^0P_{s_r(\alpha_t)}(s_r(\bi))}{Y_{\alpha_r+2\alpha_t}})\Theta_r-(\tfrac{\delta_{i_{r}}^0\delta_{i_{t}}^0(2-Y_{\alpha_r+\alpha_t})}{Y_r^2Y_{\alpha_r+2\alpha_t}}+\tfrac{\delta_{i_{r}+i_t}^0P_r(s_t(\bi))}{Y_{\alpha_r+\alpha_t}})\Theta_t\\
		&+\tfrac{\delta_{i_r}^0\delta_{i_{t}}^0P_r(\bi)}{Y_tY_{\alpha_r+\alpha_t}}+\tfrac{\delta_{i_{r}}^0\delta_{2i_t}^0P_t(\bi)}{Y_rY_{\alpha_r+2\alpha_t}}+\tfrac{\delta_{i_{r}}^0\delta_{i_{t}}^0}{Y_r^2Y_t^2}]\epsilon(\bi).
	\end{align*}
	Similarly, we have:
	\begin{align*}\Psi_t\Psi_r\Psi_t\Psi_r\epsilon(\bi)	&=[\Theta_{s_ts_rs_ts_r}-\tfrac{\delta_{i_t}^0}{Y_t}\Theta_{s_rs_ts_r}-\tfrac{\delta_{i_{r}}^0}{Y_r}\Theta_{s_ts_rs_t}+\tfrac{\delta_{i_{r}}^0\delta_{i_{t}}^0}{Y_tY_{\alpha_r+2\alpha_t}}\Theta_{s_rs_t}+\tfrac{\delta_{i_{r}}^0\delta_{i_t}^0}{Y_rY_{\alpha_r+\alpha_t}}\Theta_{s_ts_r}\\
		&-(\tfrac{\delta_{i_{r}}^0\delta_{i_t}^0}{Y_t^2Y_{\alpha_r+\alpha_t}}+\tfrac{\delta_{i_{r}+2i_t}^0P_{t}(s_r(\bi))}{Y_{\alpha_r+2\alpha_t}})\Theta_r-(\tfrac{\delta_{i_{r}}^0\delta_{i_{t}}^0(2-Y_{\alpha_r+\alpha_t})}{Y_rY^2_{\alpha_r+2\alpha_t}}+\tfrac{\delta_{i_{r}+i_t}^0P_{s_t(\alpha_r)}(s_t(\bi))}{Y_{\alpha_r+\alpha_t}})\Theta_t\\
		&+\tfrac{\delta_{i_r}^0\delta_{i_{t}}^0P_r(\bi)}{Y_tY_{\alpha_r+\alpha_t}}+\tfrac{\delta_{i_{r}}^0\delta_{2i_t}^0P_t(\bi)}{Y_rY_{\alpha_r+2\alpha_t}}+\tfrac{\delta_{i_{r}}^0\delta_{i_{t}}^0}{Y_r^2Y_t^2}]\epsilon(\bi).
	\end{align*}
	By (\ref{Theta braid relation}) we have:
	\begin{align*}
		&(\Psi_r\Psi_t\Psi_r\Psi_t-\Psi_t\Psi_r\Psi_t\Psi_r)\epsilon(\bi)\\
		&=[(\tfrac{\delta_{i_r}^0\delta_{i_{t}}^0}{Y_t^2Y_{\alpha_r+\alpha_t}}+\tfrac{\delta_{i_{r}+2i_t}^0P_t(s_r(\bi))}{Y_{\alpha_r+2\alpha_t}}-\tfrac{\delta_{i_{r}}^0\delta_{i_t}^0}{Y_tY_{\alpha_r+\alpha_t}^2}-\tfrac{\delta_{i_{r}+2i_t}^0P_{s_r(\alpha_t)}(s_r(\bi))}{Y_{\alpha_r+2\alpha_t}})\Theta_r\\
		&+(\tfrac{\delta_{i_r}^0\delta_{i_{t}}^0(2-Y_{\alpha_r+\alpha_t})}{Y_rY_{\alpha_r+2\alpha_t}^2}+\tfrac{\delta_{i_{r}+i_t}^0P_{s_t(\alpha_r)}(s_t(\bi))}{Y_{\alpha_r+\alpha_t}}-\tfrac{\delta_{i_{r}}^0\delta_{i_{t}}^0(2-Y_{\alpha_r+\alpha_t})}{Y_r^2Y_{\alpha_r+2\alpha_t}}-\tfrac{\delta_{i_{r}+i_t}^0P_r(s_t(\bi))}{Y_{\alpha_r+\alpha_t}})\Theta_t]\epsilon(\bi)\\
		&=[(\tfrac{\delta_{i_{r}}^0\delta_{i_t}^0Y_r(1-Y_t)}{Y_t^2Y_{\alpha_r+\alpha_t}^2}+\tfrac{\delta_{i_{r}+2i_t}^0(P_t(s_r(\bi))-P_{s_r(\alpha_t)}(s_r(\bi)))}{Y_{\alpha_r+2\alpha_t}})\Theta_r\\
		&+(\tfrac{\delta_{i_{r}}^0\delta_{i_{t}}^0(2-Y_{\alpha_r+\alpha_t})(Y_r-Y_{\alpha_r+2\alpha_t})}{Y_r^2Y_{\alpha_r+2\alpha_t}^2}+\tfrac{\delta_{i_r+i_{t}}^0(P_{s_t(\alpha_r)}(s_t(\bi))-P_r(s_t(\bi))}{Y_{\alpha_r+\alpha_t}})\Theta_t]\epsilon(\bi).
	\end{align*}
	There are three cases:
	
	$Case~1\colon i_r=i_t=0$. Using (\ref{P_r}), we have
	\begin{align*}&(\Psi_r\Psi_t\Psi_r\Psi_t-\Psi_t\Psi_r\Psi_t\Psi_r)\epsilon(\bi)\\
		&=[(\tfrac{Y_r(1-Y_t)}{Y_t^2Y_{\alpha_r+\alpha_t}^2}+\tfrac{Y_t^{-1}-Y_t^{-2}-(Y_{\alpha_r+\alpha_t}^{-1}-Y_{\alpha_r+\alpha_t}^{-2})}{Y_{\alpha_r+2\alpha_t}})\Theta_r\\
		&+(\tfrac{(2-Y_{\alpha_r+\alpha_t})(Y_r-Y_{\alpha_r+2\alpha_t})}{Y_r^2Y^2_{\alpha_r+2\alpha_t}}+\tfrac{Y_{\alpha_r+2\alpha_t}^{-1}-Y_{\alpha_r+2\alpha_t}^{-2}-(Y_r^{-1}-Y_{r}^{-2})}{Y_{\alpha_r+\alpha_t}})\Theta_t]\epsilon(\bi)\\
		&=[(\tfrac{Y_r-Y_rY_t}{Y_t^2Y_{\alpha_r+\alpha_t}^2}+\tfrac{(Y_{\alpha_r+\alpha_t}-Y_t)(Y_{\alpha_r+\alpha_t}Y_t-Y_{\alpha_r+\alpha_t}-Y_t)}{Y_t^2Y_{\alpha_r+\alpha_t}^2Y_{\alpha_r+2\alpha_t}})\Theta_r\\
		&+(\tfrac{(2-Y_{\alpha_r+\alpha_t})(Y_r-Y_{\alpha_r+2\alpha_t})}{Y_r^2Y^2_{\alpha_r+2\alpha_t}}+\tfrac{(Y_r-Y_{\alpha_r+2\alpha_t})(Y_rY_{\alpha_r+2\alpha_t}-Y_r-Y_{\alpha_r+2\alpha_t})}{Y_r^2Y_{\alpha_r+2\alpha_t}^2Y_{\alpha_r+\alpha_t}})\Theta_t]\epsilon(\bi)\\
		&=[(\tfrac{Y_r-Y_rY_t}{Y_t^2Y_{\alpha_r+\alpha_t}^2}-\tfrac{Y_r-Y_rY_t}{Y_t^2Y_{\alpha_r+\alpha_t}^2})\Theta_r\\
		&+(\tfrac{(2-Y_{\alpha_r+\alpha_t})(Y_r-Y_{\alpha_r+2\alpha_t})}{Y_r^2Y^2_{\alpha_r+2\alpha_t}}-\tfrac{(Y_r-Y_{\alpha_r+2\alpha_t})(2-Y_{\alpha_r+\alpha_t})}{Y_r^2Y_{\alpha_r+2\alpha_t}^2})\Theta_t]\epsilon(\bi)\\
		&=0.
	\end{align*}
	
	$Case~2\colon i_r+2i_t=0, i_r+i_t\neq 0$. Using (\ref{Psitheta}) and (\ref{P_r}) again, we have
	\begin{align*}&(\Psi_r\Psi_t\Psi_r\Psi_t-\Psi_t\Psi_r\Psi_t\Psi_r)\epsilon(\bi)\\
		&=(P_t(s_r(\bi))-P_{s_r(\alpha_t)}(s_r(\bi)))Y_{\alpha_r+2\alpha_t}^{-1}\Theta_r\epsilon(\bi)\\
		&=	\begin{cases}\tfrac{Y_t-Y_{\alpha_r+\alpha_t}(Y_{\alpha_r+\alpha_t}-1)^{-1}}{Y_{\alpha_r+2\alpha_t}}\Psi_r\epsilon(\bi) &\text{if $i_r=-2,i_t=1,e\neq 2$};\\
			\tfrac{Y_t(Y_t-1)^{-1}-Y_{\alpha_r+\alpha_t}}{Y_{\alpha_r+2\alpha_t}}\Psi_r\epsilon(\bi) &\text{if $i_r=2,i_t=-1,e\neq 2$};\\
			\tfrac{Y_t^2(Y_t-1)^{-1}-Y_{\alpha_r+\alpha_t}^2(Y_{\alpha_r+\alpha_t}-1)^{-1}}{Y_{\alpha_r+2\alpha_t}}(\Psi_r+\tfrac{1}{Y_r})\epsilon(\bi) &\text{if $i_r=0,i_t=1,e=2$};\\
			0 & \text{else}
		\end{cases}\\
		&=\begin{cases}\tfrac{1}{1-Y_{\alpha_r+\alpha_t}}\Psi_r\epsilon(\bi) &\text{if $i_r=-2,i_t=1,e\neq 2$};\\
			\tfrac{1}{Y_t-1}\Psi_r\epsilon(\bi) &\text{if $i_r=2,i_t=-1,e\neq 2$};\\
			\tfrac{1}{1-Y_{\alpha_r+\alpha_t}}(Y_r\Psi_r+1)\epsilon(\bi) &\text{if $i_r=0,i_t=1,e=2$};\\
			0 & \text{else}
		\end{cases}\\
		&=\begin{cases}(1-Y_{-\alpha_r-\alpha_t})\Psi_r\epsilon(\bi) &\text{if $i_r=-2,i_t=1,e\neq 2$};\\
			(Y_{-\alpha_t}-1)\Psi_r\epsilon(\bi) &\text{if $i_r=2,i_t=-1,e\neq 2$};\\
			(1-Y_{-\alpha_r-\alpha_t})(Y_r\Psi_r+1)\epsilon(\bi) &\text{if $i_r=0,i_t=1,e=2$};\\
			0 & \text{else}.
		\end{cases}
	\end{align*}

	$Case~3\colon i_r+i_t=0, i_r+2i_t\neq 0$. Similar to the $Case~2$, we have
	\begin{align*}&(\Psi_r\Psi_t\Psi_r\Psi_t-\Psi_t\Psi_r\Psi_t\Psi_r)\epsilon(\bi)\\
		&=(P_{s_t(\alpha_r)}(s_t(\bi))-P_r(s_t(\bi)))Y_{\alpha_r+\alpha_t}^{-1}\Theta_t\epsilon(\bi)\\
		&=	\begin{cases}\tfrac{Y_{\alpha_r+2\alpha_t}(Y_{\alpha_r+2\alpha_t}-1)^{-1}-Y_r}{Y_{\alpha_r+\alpha_t}}\Psi_t\epsilon(\bi) &\text{if $i_r=-i_t=1,e\neq 2$};\\
			\tfrac{Y_{\alpha_r+2\alpha_t}-Y_r(Y_{r}-1)^{-1}}{Y_{\alpha_r+\alpha_t}}\Psi_t\epsilon(\bi) &\text{if $i_r=-i_t=-1,e\neq 2$};\\
			\tfrac{Y_{\alpha_r+2\alpha_t}^2(Y_{\alpha_r+2\alpha_t}-1)^{-1}-Y_r^2(Y_r-1)^{-1}}{Y_{\alpha_r+\alpha_t}}\Psi_t\epsilon(\bi) &\text{if $i_r=i_t=1,e=2$};\\
			0 & \text{else}
		\end{cases}\\
		&=\begin{cases}\tfrac{2-Y_{\alpha_r+\alpha_t}}{Y_{\alpha_r+2\alpha_t}-1}\Psi_t\epsilon(\bi) &\text{if $i_r=-i_t=1,e\neq 2$};\\
			\tfrac{2-Y_{\alpha_r+\alpha_t}}{1-Y_{r}}\Psi_t\epsilon(\bi) &\text{if $i_r=-i_t=-1,e\neq 2$};\\
			\tfrac{(2-Y_{\alpha_r+\alpha_t})(Y_t-2)Y_t}{1-Y_{\alpha_r+2\alpha_t}}\Psi_t\epsilon(\bi) &\text{if $i_r=-i_t=1,e=2$};\\
			0 & \text{else}
		\end{cases}\\
		&=\begin{cases}(Y_{-\alpha_t}+Y_{-\alpha_r-2\alpha_t}-2)\Psi_t\epsilon(\bi) &\text{if $i_r=-i_t=1,e\neq 2$};\\
			(2-Y_t-Y_{-\alpha_r})\Psi_t\epsilon(\bi) &\text{if $i_r=-i_t=-1,e\neq 2$};\\
			(Y_{-\alpha_t}+Y_{-\alpha_r-2\alpha_t}-Y_{-\alpha_r}-Y_t)\Psi_t\epsilon(\bi) &\text{if $i_r=-i_t=1,e=2$};\\
			0 & \text{else}.
		\end{cases}
	\end{align*}
	Then the assertion follows.
\end{proof}

\begin{lemma}
	\label{lem:r3->t} Let $r,t\in [n]$. If $\xymatrix@C10pt{r\ar@3{-}[r]|{\rangle}\ar@3{-}[r]&t}$, then the statement (9) of Theorem \ref{thm:KLR basis} holds. \end{lemma}

\begin{proof}
	For simplicity, we denote by $\alpha_{r,1}:=s_t(\alpha_r), \alpha_{r,2}:=s_r(\alpha_{r,1}), \alpha_{r,3}:=s_t(\alpha_{r,2}), \alpha_{r,4}:=s_r(\alpha_{r,3}), \alpha_{r,5}:=s_t(\alpha_{r,4})$. We also denote by $i_{r,1}:=s_t(\bi)_r, i_{r,2}:=s_ts_r(\bi)_r, i_{r,3}:=s_ts_rs_t(\bi)_r, i_{r,4}:=s_ts_rs_ts_r(\bi)_r, i_{r,5}:=s_ts_rs_ts_rs_t(\bi)_r$. 
	Define $\alpha_{t,j}$ and $i_{t,j}$ similarly for $j=1,2,3,4,5$. Since $\xymatrix@C10pt{r\ar@3{-}[r]|{\rangle}\ar@3{-}[r]&t}$, we have 
	\begin{align*}\alpha_{r,1}&=\alpha_r+3\alpha_t=\alpha_{r,4},~~ \alpha_{r,2}=2\alpha_r+3\alpha_t= \alpha_{r,3}, ~~\alpha_{r,5}=\alpha_r;\\
		\alpha_{t,1}&=\alpha_r+\alpha_t=\alpha_{t,4},~~~~\alpha_{t,2}=\alpha_r+2\alpha_t=\alpha_{t,3},~~~~ \alpha_{t,5}=\alpha_t.\end{align*}
	Note that $i_{r,j}=i_{\alpha_{r,j}}$ and $i_{t,j}=i_{\alpha_{t,j}}$. Set $$Y_{r,j}:=Y_{\alpha_{r,j}}, \ \ Y_{t,j}:=Y_{\alpha_{t,j}}, \ \ P_{r,j}(\bi)=P_{\alpha_{r,j}}(\bi), \ \ P_{t,j}(\bi)=P_{\alpha_{t,j}}(\bi)$$
	for $j=1,2,3,4,5$.
	Then, we have
	$$Y_{r,1}=Y_{r,4}, Y_{r,2}=Y_{r,3}, Y_{r,5}=Y_r; \ \ Y_{t,1}=Y_{t,4}, Y_{t,2}=Y_{t,3}, Y_{t,5}=Y_t$$
	and $$i_{r,1}=i_{r,4}, i_{r,2}=i_{r,3}, i_{r,5}=i_r; \ \ i_{t,1}=i_{t,4}, i_{t,2}=i_{t,3}, i_{t,5}=i_t.$$
	Thus by Lemma \ref{lem:Psiepsilon}, (\ref{Thetaepsi})-(\ref{Theta^2}) and (\ref{P_r}), we have
	\begin{align*}
		&\Psi_r\Psi_t\Psi_r\Psi_t\Psi_r\Psi_t\epsilon(\bi)\\
		&=\Psi_r\epsilon(s_ts_rs_ts_rs_t(\bi))\Psi_t\epsilon(s_rs_ts_rs_t(\bi))\Psi_r\epsilon(s_ts_rs_t(\bi))\Psi_t\epsilon(s_rs_t(\bi))\Psi_r\epsilon(s_t(\bi))\Psi_t\epsilon(\bi)\\
		&=(\Theta_r-\tfrac{\delta^0_{i_{r,5}}}{Y_r})(\Theta_t-\tfrac{\delta^0_{i_{t,4}}}{Y_t})(\Theta_r-\tfrac{\delta^0_{i_{r,3}}}{Y_r})(\Theta_t-\tfrac{\delta^0_{i_{t,2}}}{Y_t})(\Theta_r-\tfrac{\delta^0_{i_{r,1}}}{Y_r})(\Theta_t-\tfrac{\delta^0_{\bi_t}}{Y_t})\epsilon(\bi)\\
		&=(\Theta_{s_rs_t}-\tfrac{\delta^0_{i_{t,4}}}{Y_{t,1}}\Theta_r-\tfrac{\delta^0_{i_{r,5}}}{Y_r}\Theta_t+\tfrac{\delta^0_{i_{r,5}}\delta^0_{i_{t,4}}}{Y_rY_t})(\Theta_{s_rs_t}-\tfrac{\delta^0_{i_{t,2}}}{Y_{t,1}}\Theta_r-\tfrac{\delta^0_{i_{r,3}}}{Y_r}\Theta_t+\tfrac{\delta^0_{i_{r,3}}\delta^0_{i_{t,2}}}{Y_rY_t})\\
		&(\Theta_{s_rs_t}-\tfrac{\delta^0_{i_t}}{Y_{t,1}}\Theta_r-\tfrac{\delta^0_{i_{r,1}}}{Y_r}\Theta_t+\tfrac{\delta^0_{i_{r,1}}\delta^0_{i_t}}{Y_rY_t})\epsilon(\bi)\\
		&=(\Theta_{s_rs_ts_rs_t}-\tfrac{\delta^0_{i_{t,2}}}{Y_{t,3}}\Theta_{s_rs_ts_r}-\tfrac{\delta^0_{i_{r,3}}P_{t,1}(s_t(\bi))}{Y_{r,2}}\Theta_r+\tfrac{\delta^0_{i_{r,3}}\delta^0_{i_{t,2}}}{Y_{r,2}Y_{s_rs_t(\alpha_t)}}\Theta_{s_rs_t}-\tfrac{\delta^0_{i_{t,4}}P_r(s_ts_rs_t(\bi))}{Y_{t,1}}\Theta_t\\
		&+\tfrac{\delta^0_{i_{t,4}}\delta^0_{i_{t,2}}P_r(s_rs_t(\bi))}{Y_{t}Y_{t,1}}+\tfrac{\delta^0_{i_{t,4}}\delta^0_{i_{r,3}}}{Y_{s_r(\alpha_r)}Y_{t,1}}\Theta_{s_rs_t}-\tfrac{\delta^0_{i_{t,4}}\delta_{i_{r,3}}^0\delta^0_{i_{t,2}}}{Y_{s_r(\alpha_r)}Y_{t,1}^2}\Theta_r-\tfrac{\delta^0_{i_{r,5}}}{Y_r}\Theta_{s_ts_rs_t}+\tfrac{\delta^0_{i_{r,5}}\delta^0_{i_{t,2}}}{Y_rY_{t,2}}\Theta_{s_ts_r}\\
		&+\tfrac{\delta^0_{i_{r,5}}\delta^0_{i_{r,3}}P_t(s_rs_t(\bi))}{Y_rY_{r,1}}
		-\tfrac{\delta^0_{i_{r,5}}\delta_{i_{r,3}}^0\delta^0_{i_{t,2}}}{Y_rY_{r,1}Y_{s_t(\alpha_t)}}\Theta_t+\tfrac{\delta^0_{i_{r,5}}\delta^0_{i_{t,4}}}{Y_rY_t}\Theta_{s_rs_t}-\tfrac{\delta^0_{i_{r,5}}\delta^0_{i_{t,4}}\delta^0_{i_{t,2}}}{Y_rY_tY_{t,1}}\Theta_r-\tfrac{\delta^0_{i_{r,5}}\delta^0_{i_{t,4}}\delta^0_{i_{r,3}}}{Y_r^2Y_t}\Theta_t\\
		&+\tfrac{\delta^0_{i_{r,5}}\delta^0_{i_{t,4}}\delta^0_{i_{r,3}}\delta^0_{i_{t,2}}}{Y_r^2Y_t^2})(\Theta_{s_rs_t}-\tfrac{\delta^0_{i_t}}{Y_{t,1}}\Theta_r-\tfrac{\delta^0_{i_{r,1}}}{Y_r}\Theta_t+\tfrac{\delta^0_{i_{r,1}}\delta^0_{i_t}}{Y_rY_t})\epsilon(\bi)\\
		&=[\Theta_{s_rs_ts_rs_t}-\tfrac{\delta^0_{i_{t,2}}}{Y_{t,2}}\Theta_{s_rs_ts_r}-\tfrac{\delta^0_{i_r}}{Y_r}\Theta_{s_ts_rs_t}+\tfrac{\delta^0_{i_r}\delta^0_{i_t}(2-Y_{t,2})}{Y_tY_{r,2}}\Theta_{s_rs_t}+\tfrac{\delta^0_{i_r}\delta^0_{i_{t,2}}}{Y_rY_{t,2}}\Theta_{s_ts_r}\\
		&-(\tfrac{\delta^0_{i_{r,2}}P_{t,1}(s_t(\bi))}{Y_{r,2}}+\tfrac{\delta^0_{i_r}\delta^0_{i_t}}{Y_tY_{t,1}^2})\Theta_r-(\tfrac{\delta^0_{i_{t,1}}P_r(s_ts_rs_t(\bi))}{Y_{t,1}}+\tfrac{\delta^0_{i_r}\delta^0_{i_t}(3-Y_{t,2}-Y_{t,1})}{Y_r^2Y_{r,1}})\Theta_t+\tfrac{\delta^0_{i_r}\delta^0_{i_t}}{Y_rY_t^2Y_{t,1}}\\
		&+\tfrac{\delta^0_{i_r}\delta^0_{i_{r,1}}P_t(s_rs_t(\bi))}{Y_rY_{r,1}}](\Theta_{s_rs_t}-\tfrac{\delta^0_{i_t}}{Y_{t,1}}\Theta_r-\tfrac{\delta^0_{i_{r,1}}}{Y_r}\Theta_t+\tfrac{\delta^0_{i_{r}}\delta^0_{i_t}}{Y_rY_t})\epsilon(\bi)\\
		&=[\Theta_{s_rs_ts_rs_ts_rs_t}-\tfrac{\delta^0_{i_t}}{Y_{t,5}}\Theta_{s_rs_ts_rs_ts_r}-\tfrac{\delta^0_{i_{r,1}}P_{t,3}(s_rs_ts_r(\bi))}{Y_{r,4}}\Theta_{s_rs_ts_r}+\tfrac{\delta^0_{i_r}\delta^0_{i_t}}{Y_{r,4}Y_{s_rs_ts_rs_t(\alpha_t)}}\Theta_{s_rs_ts_rs_t}\\
		&-\tfrac{\delta^0_{i_{t,2}}P_{r,2}(s_r(\bi))P_{t,1}(s_r(\bi))}{Y_{t,2}}\Theta_r+\tfrac{\delta^0_{i_{t,2}}\delta^0_{i_t}P_{r,2}(s_rs_t(\bi))}{Y_{t,2}Y_{s_rs_t(\alpha_t)}}\Theta_{s_rs_t}+\tfrac{\delta^0_{i_{r,1}}\delta^0_{i_{t,2}}}{Y_{t,2}Y_{s_rs_ts_r(\alpha_r)}}\Theta_{s_rs_ts_rs_t}\\
		&-\tfrac{\delta^0_{i_r}\delta^0_{i_t}}{Y_{s_rs_ts_r(\alpha_r)}Y_{t,2}Y_{t,3}}\Theta_{s_rs_ts_r}-\tfrac{\delta^0_{i_r}}{Y_r}\Theta_{s_ts_rs_ts_rs_t}+\tfrac{\delta^0_{i_r}\delta^0_{i_t}}{Y_rY_{t,4}}\Theta_{s_ts_rs_ts_r}+\tfrac{\delta^0_{i_r}\delta^0_{i_{r,2}}P_{t,2}(s_ts_r(\bi))}{Y_rY_{r,3}}\Theta_{s_ts_r}\\
		&-\tfrac{\delta^0_{i_r}\delta^0_{i_t}}{Y_rY_{r,3}Y_{s_ts_rs_t(\alpha_t)}}\Theta_{s_ts_rs_t}+\tfrac{\delta^0_{i_r}\delta^0_{i_t}(2-Y_{t,2})}{Y_tY_{r,2}}\Theta_{s_rs_ts_rs_t}-\tfrac{\delta^0_{i_r}\delta^0_{i_t}(2-Y_{t,2})}{Y_tY_{r,2}Y_{t,3}}\Theta_{s_rs_ts_r}\\
		&-\tfrac{\delta^0_{i_r}\delta^0_{i_t}(2-Y_{t,2})P_{t,1}(s_r(\bi))}{Y_tY_{r,2}^2}\Theta_r+\tfrac{\delta^0_{i_r}\delta^0_{i_t}(2-Y_{t,2})}{Y_tY_{r,2}^2Y_{s_rs_t(\alpha_t)}}\Theta_{s_rs_t}+\tfrac{\delta^0_{i_r}\delta^0_{i_{t,2}}P_{r,1}(\bi)P_t(\bi)}{Y_rY_{t,2}}-\tfrac{\delta^0_{i_r}\delta^0_{i_t}P_{r,1}(s_t(\bi))}{Y_rY_{t,2}Y_{s_t(\alpha_t)}}\Theta_t\\
		&-\tfrac{\delta^0_{i_r}\delta^0_{i_t}\delta_{i_{r,1}}^0}{Y_rY_{t,2}Y_{s_ts_r(\alpha_r)}}\Theta_{s_ts_rs_t}+\tfrac{\delta^0_{i_r}\delta^0_{i_t}}{Y_rY_{t,2}^2Y_{s_ts_r(\alpha_r)}}\Theta_{s_ts_r}-(\tfrac{\delta^0_{i_{r,2}}P_{t,1}(s_t(\bi))}{Y_{r,2}}+\tfrac{\delta^0_{i_r}\delta^0_{i_t}}{Y_tY_{t,1}^2})P_{r}(s_t(\bi))\Theta_t\\
		&+(\tfrac{\delta^0_{i_{r,2}}\delta^0_{i_t}P_{t,1}(s_t(\bi))}{Y_{r,2}Y_t}+\tfrac{\delta^0_{i_r}\delta^0_{i_t}}{Y_t^2Y_{t,1}^2})P_{r}(\bi)+(\tfrac{\delta^0_{i_{r,1}}\delta^0_{i_{r,2}}P_{t,1}(s_t(\bi))}{Y_{s_r(\alpha_r)}Y_{r,2}}+\tfrac{\delta^0_{i_r}\delta^0_{i_t}}{Y_{s_r(\alpha_r)}Y_tY_{t,1}^2})\Theta_{s_rs_t}\\
		&-(\tfrac{\delta^0_{i_r}\delta^0_{i_t}P_{t,1}(s_t(\bi))}{Y_{s_r(\alpha_r)}Y_{r,2}Y_{t,1}}+\tfrac{\delta^0_{i_r}\delta^0_{i_t}}{Y_{s_r(\alpha_r)}Y_tY_{t,1}^3})\Theta_r-(\tfrac{\delta^0_{i_{t,1}}P_r(s_ts_rs_t(\bi))}{Y_{t,1}}+\tfrac{\delta^0_{i_r}\delta^0_{i_t}(3-Y_{t,2}-Y_{t,1})}{Y_r^2Y_{r,1}})\Theta_{s_ts_rs_t}\\
		&+(\tfrac{\delta^0_{i_{t,1}}\delta^0_{i_t}P_r(s_ts_rs_t(\bi))}{Y_{t,1}Y_{t,2}}+\tfrac{\delta^0_{i_r}\delta^0_{i_t}(3-Y_{t,2}-Y_{t,1})}{Y_r^2Y_{r,1}Y_{t,2}})\Theta_{s_ts_r}+\tfrac{\delta^0_{i_{r,1}}\delta^0_{i_{t,1}}P_r(s_ts_rs_t(\bi))P_t(\bi)}{Y_{r,1}Y_{t,1}}\\
		&+\tfrac{\delta^0_{i_r}\delta^0_{i_t}(3-Y_{t,2}-Y_{t,1})P_t(\bi)}{Y_r^2Y_{r,1}^2}-(\tfrac{\delta^0_{i_r}\delta^0_{i_t}P_r(s_ts_rs_t(\bi))}{Y_{t,1}Y_{r,1}Y_{s_t(\alpha_t)}}+\tfrac{\delta^0_{i_r}\delta^0_{i_t}(3-Y_{t,2}-Y_{t,1})}{Y_r^2Y_{r,1}^2Y_{s_t(\alpha_t)}})\Theta_t\\
		&+(\tfrac{\delta^0_{i_r}\delta^0_{i_t}}{Y_rY_t^2Y_{t,1}}+\tfrac{\delta^0_{i_r}\delta^0_{i_{r,1}}P_t(s_rs_t(\bi))}{Y_rY_{r,1}})\Theta_{s_rs_t}-(\tfrac{\delta^0_{i_r}\delta^0_{i_t}}{Y_rY_t^2Y_{t,1}^2}+\tfrac{\delta^0_{i_r}\delta^0_{i_t}P_t(s_rs_t(\bi))}{Y_rY_{r,1}Y_{t,1}})\Theta_r\\
		&-(\tfrac{\delta^0_{i_r}\delta^0_{i_t}}{Y_r^2Y_t^2Y_{t,1}}+\tfrac{\delta^0_{i_r}\delta^0_{i_{r,1}}P_t(s_rs_t(\bi))}{Y_r^2Y_{r,1}})\Theta_t+\tfrac{\delta^0_{i_r}\delta^0_{i_t}}{Y_r^2Y_t^3Y_{t,1}}+\tfrac{\delta^0_{i_r}\delta^0_{i_{t}}P_t(s_rs_t(\bi))}{Y_r^2Y_tY_{r,1}}]\epsilon(\bi)\\
		&=[\Theta_{s_rs_ts_rs_ts_rs_t}-\tfrac{\delta_{i_t}^0}{Y_t}\Theta_{s_rs_ts_rs_ts_r}-\tfrac{\delta_{i_r}^0}{Y_r}\Theta_{s_ts_rs_ts_rs_t}+\tfrac{\delta_{i_{r}}^0\delta_{i_{t}}^0}{Y_tY_{r,1}}\Theta_{s_rs_ts_rs_t}+\tfrac{\delta_{i_r}^0\delta_{i_t}^0}{Y_rY_{t,1}}\Theta_{s_ts_rs_ts_r}\\
		&-(\tfrac{\delta_{i_r}^0\delta_{i_{t}}^0}{Y_tY_{t,2}^2}+\tfrac{\delta_{i_{r,1}}^0P_{t,2}(s_rs_ts_r(\bi))}{Y_{r,1}})\Theta_{s_rs_ts_r}-(\tfrac{\delta_{i_r}^0\delta_{i_{t}}^0(3-3Y_{t,1}+Y_{t,1}^2)}{Y_r^2Y_{r,2}}+\tfrac{\delta^0_{i_{t,1}}P_r(s_ts_rs_t(\bi))}{Y_{t,1}})\Theta_{s_ts_rs_t}\\			
		&+(\tfrac{\delta_{i_r}^0\delta_{i_{t}}^0[Y_{t,1}^2+Y_{t,2}(Y_{r,2}-Y_{t,1})]}{Y_{r,2}Y_t^2Y_{t,1}^2Y_{t,2}}+\tfrac{\delta_{i_r}^0\delta_{i_{r,1}}^0P_{t}(s_rs_t(\bi))}{Y_rY_{r,1}}+\tfrac{\delta_{i_r}^0\delta_{i_{r,2}}^0P_{t,1}(s_t(\bi))}{Y_{-\alpha_r}Y_{r,2}})\Theta_{s_rs_t}+(\tfrac{\delta_{i_r}^0\delta_{i_t}^0}{Y_rY_{t,1}Y_{t,2}^2}\\
		&+\tfrac{\delta_{i_r}^0\delta_{i_{r,2}}^0P_{t,2}(s_ts_r(\bi))}{Y_rY_{r,2}})\Theta_{s_ts_r}-(\tfrac{\delta_{i_{r}}^0\delta_{i_{t}}^0(2-Y_{{t,2}})[(Y_{t,1}-1)Y_{r,1}+Y_{r,2}]}{Y_{r,1}Y_{r,2}^2Y_tY_{t,1}^2}+\tfrac{\delta_{i_{t,2}}^0P_{r,2}(s_r(\bi))P_{t,1}(s_r(\bi))}{Y_{t,2}})\Theta_r\\
		&-(\tfrac{\delta_{i_r}^0\delta_{i_{t}}^0[Y_t^2+Y_{t,2}(Y_{r,1}-Y_t)]}{Y_rY_{r,1}Y_t^2Y_{t,1}^2Y_{t,2}}+\tfrac{\delta_{i_{r,2}}^0P_r(s_t(\bi))P_{t,1}(s_t(\bi))}{Y_{r,2}}+\tfrac{\delta_{i_r}^0\delta_{i_{r,1}}^0P_t(s_rs_t(\bi))}{Y_r^2Y_{r,1}})\Theta_t
		+\tfrac{\delta_{i_{r,2}}^0\delta_{i_t}^0P_r(\bi)P_{t,1}(s_t(\bi))}{Y_tY_{r,2}}\\
		&+\tfrac{\delta_{i_r}^0\delta_{i_{t,2}}^0P_{r,1}(\bi)P_t(\bi)}{Y_rY_{t,2}}+\tfrac{\delta_{i_{r,1}}^0\delta_{t,1}^0P_{r}(s_ts_rs_t(\bi))P_t(\bi)}{Y_{r,1}Y_{t,1}}+\tfrac{\delta_{i_r}^0\delta_{i_t}^0(Y_r-1)}{Y_r^2Y_{t}^2Y_{t,1}^2}+\tfrac{\delta_{i_r}^0\delta_{i_t}^0(Y_t-1)(3-Y_{t,1}-Y_{t,2})}{Y_r^2Y_{r,1}^2Y_{t}^2}\\
		&+\tfrac{\delta_{i_r}^0\delta_{i_t}^0}{Y_r^2Y_{t}^3Y_{t,1}}+\tfrac{\delta_{i_r}^0\delta_{i_t}^0(Y_t-1)}{Y_r^2Y_{t}^3Y_{r,1}}]\epsilon(\bi)
	\end{align*}	
	Similarly, we have
	\begin{align*}
		&\Psi_t\Psi_r\Psi_t\Psi_r\Psi_t\Psi_r\epsilon(\bi)\\
		&=\Psi_t\epsilon(s_rs_ts_rs_ts_r(\bi))\Psi_r\epsilon(s_ts_rs_ts_r(\bi))\Psi_t\epsilon(s_rs_ts_r(\bi))\Psi_r\epsilon(s_ts_r(\bi))\Psi_t\epsilon(s_r(\bi))\Psi_r\epsilon(\bi)\\
		&=(\Theta_t-\tfrac{\delta^0_{i_{t,5}}}{Y_t})(\Theta_r-\tfrac{\delta^0_{i_{r,4}}}{Y_r})(\Theta_t-\tfrac{\delta^0_{i_{t,3}}}{Y_t})(\Theta_r-\tfrac{\delta^0_{i_{r,2}}}{Y_r})(\Theta_t-\tfrac{\delta^0_{i_{t,1}}}{Y_t})(\Theta_r-\tfrac{\delta^0_{\bi_r}}{Y_r})\epsilon(\bi)\\
		&=(\Theta_{s_ts_r}-\tfrac{\delta^0_{i_{r,4}}}{Y_{r,1}}\Theta_t-\tfrac{\delta^0_{i_{t,5}}}{Y_t}\Theta_r+\tfrac{\delta^0_{i_{t,5}}\delta^0_{i_{r,4}}}{Y_rY_t})(\Theta_{s_ts_r}-\tfrac{\delta^0_{i_{r,2}}}{Y_{r,1}}\Theta_t-\tfrac{\delta^0_{i_{t,3}}}{Y_t}\Theta_r+\tfrac{\delta^0_{i_{t,3}}\delta^0_{i_{r,2}}}{Y_rY_t})\\
		&(\Theta_{s_ts_r}-\tfrac{\delta^0_{i_r}}{Y_{r,1}}\Theta_t-\tfrac{\delta^0_{i_{t,1}}}{Y_t}\Theta_r+\tfrac{\delta^0_{i_{t,1}}\delta^0_{i_r}}{Y_rY_t})\epsilon(\bi)\\
		&=(\Theta_{s_ts_rs_ts_r}-\tfrac{\delta^0_{i_{r,2}}}{Y_{r,3}}\Theta_{s_ts_rs_t}-\tfrac{\delta^0_{i_{t,3}}P_{r,1}(s_r(\bi))}{Y_{t,2}}\Theta_t+\tfrac{\delta^0_{i_{t,3}}\delta^0_{i_{r,2}}}{Y_{s_ts_r(\alpha_r)}Y_{t,2}}\Theta_{s_ts_r}-\tfrac{\delta^0_{i_{r,4}}P_t(s_rs_ts_r(\bi))}{Y_{r,1}}\Theta_r\\
		&+\tfrac{\delta^0_{i_{r,4}}\delta^0_{i_{r,2}}P_t(s_ts_r(\bi))}{Y_{r}Y_{r,1}}+\tfrac{\delta^0_{i_{r,4}}\delta^0_{i_{t,3}}}{Y_{r,1}Y_{s_t(\alpha_t)}}\Theta_{s_ts_r}-\tfrac{\delta^0_{i_{r,4}}\delta_{i_{t,3}}^0\delta^0_{i_{r,2}}}{Y_{r,1}^2Y_{s_t(\alpha_t)}}\Theta_t-\tfrac{\delta^0_{i_{t,5}}}{Y_t}\Theta_{s_rs_ts_r}+\tfrac{\delta^0_{i_{t,5}}\delta^0_{i_{r,2}}}{Y_{r,2}Y_t}\Theta_{s_rs_t}\\
		&+\tfrac{\delta^0_{i_{t,5}}\delta^0_{i_{t,3}}P_r(s_ts_r(\bi))}{Y_tY_{t,1}}
		-\tfrac{\delta^0_{i_{t,5}}\delta_{i_{t,3}}^0\delta^0_{i_{r,2}}}{Y_{s_r(\alpha_r)}Y_tY_{t,1}}\Theta_r+\tfrac{\delta^0_{i_{t,5}}\delta^0_{i_{r,4}}}{Y_rY_t}\Theta_{s_ts_r}-\tfrac{\delta^0_{i_{t,5}}\delta^0_{i_{r,4}}\delta^0_{i_{r,2}}}{Y_rY_{r,1}Y_t}\Theta_t-\tfrac{\delta^0_{i_{t,5}}\delta^0_{i_{r,4}}\delta^0_{i_{t,3}}}{Y_rY_t^2}\Theta_r\\
		&+\tfrac{\delta^0_{i_{t,5}}\delta^0_{i_{r,4}}\delta^0_{i_{t,3}}\delta^0_{i_{r,2}}}{Y_r^2Y_t^2})(\Theta_{s_ts_r}-\tfrac{\delta^0_{i_r}}{Y_{r,1}}\Theta_t-\tfrac{\delta^0_{i_{t,1}}}{Y_t}\Theta_r+\tfrac{\delta^0_{i_{t,1}}\delta^0_{i_r}}{Y_rY_t})\epsilon(\bi)\\
		&=[\Theta_{s_ts_rs_ts_r}-\tfrac{\delta^0_{i_t}}{Y_t}\Theta_{s_rs_ts_r}-\tfrac{\delta^0_{i_{r,2}}}{Y_{r,2}}\Theta_{s_ts_rs_t}+\tfrac{\delta^0_{i_t}\delta^0_{i_{r,2}}}{Y_tY_{r,2}}\Theta_{s_rs_t}+\tfrac{\delta^0_{i_r}\delta^0_{i_t}(2-Y_{t,1})}{Y_rY_{t,2}}\Theta_{s_ts_r}\\
		&-(\tfrac{\delta^0_{i_{r,1}}P_{t}(s_rs_ts_r(\bi))}{Y_{r,1}}+\tfrac{\delta^0_{i_r}\delta^0_{i_t}}{Y_t^2Y_{t,1}})\Theta_r-(\tfrac{\delta^0_{i_{t,2}}P_{r,1}(s_r(\bi))}{Y_{t,2}}+\tfrac{\delta^0_{i_r}\delta^0_{i_t}(3-Y_{t,2}-Y_{t,1})}{Y_rY_{r,1}^2})\Theta_t+\tfrac{\delta^0_{i_r}\delta^0_{i_t}}{Y_rY_t^2Y_{t,1}}\\
		&+\tfrac{\delta^0_{i_{r}}\delta^0_{i_{r,1}}P_t(s_ts_r(\bi))}{Y_rY_{r,1}}](\Theta_{s_ts_r}-\tfrac{\delta^0_{i_r}}{Y_{r,1}}\Theta_t-\tfrac{\delta^0_{i_{t,1}}}{Y_t}\Theta_r+\tfrac{\delta^0_{i_{r}}\delta^0_{i_t}}{Y_rY_t})\epsilon(\bi)\\
		&=[\Theta_{s_ts_rs_ts_rs_ts_r}-\tfrac{\delta^0_{i_r}}{Y_{r,5}}\Theta_{s_ts_rs_ts_rs_t}-\tfrac{\delta^0_{i_{t,1}}P_{r,3}(s_ts_rs_t(\bi))}{Y_{t,4}}\Theta_{s_ts_rs_t}+\tfrac{\delta^0_{i_r}\delta^0_{i_t}}{Y_{s_ts_rs_ts_r(\alpha_r)}Y_{t,4}}\Theta_{s_ts_rs_ts_r}\\
		&-\tfrac{\delta^0_{i_t}}{Y_t}\Theta_{s_rs_ts_rs_ts_r}+\tfrac{\delta^0_{i_r}\delta^0_{i_t}}{Y_{r,4}Y_t}\Theta_{s_rs_ts_rs_t}+\tfrac{\delta^0_{i_t}\delta^0_{i_{t,1}}P_{r,2}(s_rs_t(\bi))}{Y_tY_{t,3}}\Theta_{s_rs_t}-\tfrac{\delta^0_{i_r}\delta^0_{i_t}}{Y_{s_rs_ts_r(\alpha_r)}Y_{t}Y_{t,3}}\Theta_{s_rs_ts_r}\\
		&-\tfrac{\delta^0_{i_{r,2}}P_{t,2}(s_t(\bi))P_{r,1}(s_t(\bi))}{Y_{r,2}}\Theta_t+\tfrac{\delta^0_{i_{r,2}}\delta^0_{i_r}P_{t,2}(s_ts_r(\bi))}{Y_{r,2}Y_{s_ts_r(\alpha_r)}}\Theta_{s_ts_r}+\tfrac{\delta^0_{i_{t,1}}\delta^0_{i_{r,2}}}{Y_{r,2}Y_{s_ts_rs_t(\alpha_t)}}\Theta_{s_ts_rs_ts_r}-\tfrac{\delta^0_{i_r}\delta^0_{i_t}}{Y_{r,2}Y_{r,3}Y_{s_ts_rs_t(\alpha_t)}}\Theta_{s_ts_rs_t}\\
		&+\tfrac{\delta^0_{i_t}\delta^0_{i_{r,2}}P_{t,1}(\bi)P_r(\bi)}{Y_tY_{r,2}}-\tfrac{\delta^0_{i_r}\delta^0_{i_t}P_{t,1}(s_r(\bi))}{Y_tY_{r,2}Y_{s_r(\alpha_r)}}\Theta_r-\tfrac{\delta^0_{i_t}\delta^0_{i_{r,2}}\delta_{i_{t,1}}^0}{Y_tY_{r,2}Y_{s_rs_t(\alpha_t)}}\Theta_{s_rs_ts_r}+\tfrac{\delta^0_{i_r}\delta^0_{i_t}}{Y_tY_{r,2}^2Y_{s_rs_t(\alpha_t)}}\Theta_{s_rs_t}\\
		&+\tfrac{\delta^0_{i_r}\delta^0_{i_t}(2-Y_{t,1})}{Y_rY_{t,2}}\Theta_{s_ts_rs_ts_r}-\tfrac{\delta^0_{i_r}\delta^0_{i_t}(2-Y_{t,1})}{Y_rY_{t,2}Y_{r,3}}\Theta_{s_ts_rs_t}-\tfrac{\delta^0_{i_r}\delta^0_{i_t}(2-Y_{t,1})P_{r,1}(s_t(\bi))}{Y_rY_{t,2}^2}\Theta_t+\tfrac{\delta^0_{i_r}\delta^0_{i_t}(2-Y_{t,1})}{Y_rY_{t,2}^2Y_{s_ts_r(\alpha_r)}}\Theta_{s_ts_r}\\
		&-(\tfrac{\delta^0_{i_{r,1}}P_t(s_rs_ts_r(\bi))}{Y_{r,1}}+\tfrac{\delta^0_{i_r}\delta^0_{i_t}}{Y_t^2Y_{t,1}})\Theta_{s_rs_ts_r}+(\tfrac{\delta^0_{i_{r,1}}\delta^0_{i_r}P_t(s_rs_ts_r(\bi))}{Y_{r,1}Y_{r,2}}+\tfrac{\delta^0_{i_r}\delta^0_{i_t}}{Y_t^2Y_{t,1}Y_{r,2}})\Theta_{s_rs_t}+\tfrac{\delta^0_{i_{r,1}}\delta^0_{i_{t,1}}P_t(s_rs_ts_r(\bi))P_r(\bi)}{Y_{r,1}Y_{t,1}}\\
		&+\tfrac{\delta^0_{i_r}\delta^0_{i_t}P_r(\bi)}{Y_t^2Y_{t,1}^2}-(\tfrac{\delta^0_{i_r}\delta^0_{i_t}P_t(s_rs_ts_r(\bi))}{Y_{t,1}Y_{r,1}Y_{s_r(\alpha_r)}}+\tfrac{\delta^0_{i_r}\delta^0_{i_t}}{Y_t^2Y_{t,1}^2Y_{s_r(\alpha_r)}})\Theta_r-(\tfrac{\delta^0_{i_{t,2}}P_{r,1}(s_r(\bi))}{Y_{t,2}}+\tfrac{\delta^0_{i_r}\delta^0_{i_t}(3-Y_{t,2}-Y_{t,1})}{Y_rY_{r,1}^2})P_{t}(s_r(\bi))\Theta_r\\
		&+(\tfrac{\delta^0_{i_{t,2}}\delta^0_{i_r}P_{r,1}(s_r(\bi))}{Y_{t,2}Y_r}+\tfrac{\delta^0_{i_r}\delta^0_{i_t}(3-Y_{t,2}-Y_{t,1})}{Y_r^2Y_{r,1}^2})P_{t}(\bi)+(\tfrac{\delta^0_{i_{t,1}}\delta^0_{i_{t,2}}P_{r,1}(s_r(\bi))}{Y_{t,2}Y_{s_t(\alpha_t)}}+\tfrac{\delta^0_{i_r}\delta^0_{i_t}(3-Y_{t,2}-Y_{t,1})}{Y_rY_{r,1}^2Y_{s_t(\alpha_t)}})\Theta_{s_ts_r}\\
		&-(\tfrac{\delta^0_{i_r}\delta^0_{i_t}P_{r,1}(s_r(\bi))}{Y_{r,1}Y_{t,2}Y_{s_t(\alpha_t)}}+\tfrac{\delta^0_{i_r}\delta^0_{i_t}(3-Y_{t,2}-Y_{t,1})}{Y_rY_{r,1}^3Y_{s_t(\alpha_t)}})\Theta_t\\
		&+(\tfrac{\delta^0_{i_r}\delta^0_{i_t}}{Y_rY_t^2Y_{t,1}}+\tfrac{\delta^0_{i_r}\delta^0_{i_{r,1}}P_t(s_ts_r(\bi))}{Y_rY_{r,1}})\Theta_{s_ts_r}-(\tfrac{\delta^0_{i_r}\delta^0_{i_t}}{Y_rY_t^2Y_{t,1}Y_{r,1}}+\tfrac{\delta^0_{i_r}\delta^0_{i_{r,1}}P_t(s_ts_r(\bi))}{Y_rY_{r,1}^2})\Theta_t\\
		&-(\tfrac{\delta^0_{i_r}\delta^0_{i_t}}{Y_rY_t^3Y_{t,1}}+\tfrac{\delta^0_{i_r}\delta^0_{i_{t}}P_t(s_ts_r(\bi))}{Y_rY_{r,1}Y_{t}})\Theta_r+\tfrac{\delta^0_{i_r}\delta^0_{i_t}}{Y_r^2Y_t^3Y_{t,1}}+\tfrac{\delta^0_{i_r}\delta^0_{i_{t}}P_t(s_ts_r(\bi))}{Y_r^2Y_tY_{r,1}}]\epsilon(\bi)\\
		&=[\Theta_{s_ts_rs_ts_rs_ts_r}-\tfrac{\delta_{i_t}^0}{Y_t}\Theta_{s_rs_ts_rs_ts_r}-\tfrac{\delta_{i_r}^0}{Y_r}\Theta_{s_ts_rs_ts_rs_t}+\tfrac{\delta_{i_{r}}^0\delta_{i_{t}}^0}{Y_tY_{r,1}}\Theta_{s_rs_ts_rs_t}+\tfrac{\delta_{i_r}^0\delta_{i_t}^0}{Y_rY_{t,1}}\Theta_{s_ts_rs_ts_r}\\
		&-(\tfrac{\delta_{i_r}^0\delta_{i_{t}}^0}{Y_t^2Y_{t,2}}+\tfrac{\delta_{i_{r,1}}^0P_{t}(s_rs_ts_r(\bi))}{Y_{r,1}})\Theta_{s_rs_ts_r}-(\tfrac{\delta_{i_r}^0\delta_{i_{t}}^0(3-3Y_{t,1}+Y_{t,1}^2)}{Y_rY_{r,2}^2}+\tfrac{\delta^0_{i_{t,1}}P_{r,2}(s_ts_rs_t(\bi))}{Y_{t,1}})\Theta_{s_ts_rs_t}\\			
		&+(\tfrac{\delta_{i_r}^0\delta_{i_t}^0}{Y_{r,2}Y_{t,2}Y_{t}^2}+\tfrac{\delta_{i_r}^0\delta_{i_{r,1}}^0P_{t}(s_rs_t(\bi))}{Y_{r,1}Y_{r,2}})\Theta_{s_rs_t}+(\tfrac{\delta_{i_r}^0\delta_{i_{t}}^0[Y_t^2+Y_{t,1}(Y_{t,2}-Y_t)]}{Y_{r}Y_t^2Y_{t,1}Y_{t,2}^2}+\tfrac{\delta_{i_r}^0\delta_{i_{r,1}}^0P_{t}(s_ts_r(\bi))}{Y_rY_{r,1}}\\
		&+\tfrac{\delta_{i_r}^0\delta_{i_{r,2}}^0P_{t,2}(s_ts_r(\bi))}{Y_{-\alpha_{r,1}}Y_{r,2}})\Theta_{s_ts_r}
		-(\tfrac{\delta_{i_{r}}^0\delta_{i_{t}}^0(2-Y_{{t,2}})[(Y_{r,1}-1)Y_{t,1}+Y_{t}]}{Y_{r,1}^2Y_{r,2}Y_t^2Y_{t,1}}+\tfrac{\delta_{i_{t,2}}^0P_{r,1}(s_r(\bi))P_{t}(s_r(\bi))}{Y_{t,2}})\Theta_r\\
		&-(\tfrac{\delta_{i_r}^0\delta_{i_{t}}^0[Y_t^2+Y_{t,1}(Y_{t,2}-Y_t)]}{Y_rY_{r,1}Y_t^2Y_{t,1}Y_{t,2}^2}+\tfrac{\delta_{i_{r,2}}^0P_{t,2}(s_t(\bi))P_{r,1}(s_t(\bi))}{Y_{r,2}}+\tfrac{\delta_{i_r}^0\delta_{i_{r,1}}^0P_t(s_ts_r(\bi))}{Y_rY_{r,1}^2})\Theta_t
		+\tfrac{\delta_{i_{r,2}}^0\delta_{i_t}^0P_r(\bi)P_{t,1}(\bi)}{Y_tY_{r,2}}\\
		&+\tfrac{\delta_{i_r}^0\delta_{i_{t,2}}^0P_{r,1}(s_r(\bi))P_t(\bi)}{Y_rY_{t,2}}+\tfrac{\delta_{i_{r,1}}^0\delta_{t,1}^0P_{t}(s_rs_ts_r(\bi))P_r(\bi)}{Y_{r,1}Y_{t,1}}+\tfrac{\delta_{i_r}^0\delta_{i_t}^0(Y_r-1)}{Y_r^2Y_{t}^2Y_{t,1}^2}+\tfrac{\delta_{i_r}^0\delta_{i_t}^0(Y_t-1)(3-Y_{t,1}-Y_{t,2})}{Y_r^2Y_{r,1}^2Y_{t}^2}\\
		&+\tfrac{\delta_{i_r}^0\delta_{i_t}^0}{Y_r^2Y_{t}^3Y_{t,1}}+\tfrac{\delta_{i_r}^0\delta_{i_t}^0(Y_t-1)}{Y_r^2Y_{t}^3Y_{r,1}}]\epsilon(\bi)
	\end{align*}

	Using relations (\ref{Theta braid relation}) and (\ref{P_r}), we have
	\begin{align*}
		&(\Psi_r\Psi_t\Psi_r\Psi_t\Psi_r\Psi_t-\Psi_t\Psi_r\Psi_t\Psi_r\Psi_t\Psi_r)\epsilon(\bi)\\
		&=[(\tfrac{\delta_{i_{r}}^0\delta_{i_t}^0(Y_{t,2}-Y_{t})}{Y_t^2Y_{t,2}^2}+\tfrac{\delta_{i_{r,1}}^0(P_t(s_rs_ts_r(\bi))-P_{t,2}(s_rs_ts_r(\bi)))}{Y_{r,1}})\Theta_{s_rs_ts_r}\\
		&+(\tfrac{\delta_{i_{r}}^0\delta_{i_t}^0(Y_r-Y_{r,2})(3-3Y_{t,1}+Y_{t,1}^2)}{Y_r^2Y_{r,2}^2}+\tfrac{\delta_{i_{t,1}}^0(P_{r,2}(s_ts_rs_t(\bi))-P_{r}(s_ts_rs_t(\bi)))}{Y_{t,1}})\Theta_{s_ts_rs_t}\\
		&+(\tfrac{\delta_{i_r}^0\delta_{i_t}^0(Y_{r,2}-Y_{t,1})}{Y_{r,2}Y_t^2Y_{t,1}^2}+\tfrac{\delta_{i_r}^0\delta_{i_{r,1}}^0(P_{t,1}(s_t(\bi))-P_{t}(s_rs_t(\bi)))}{Y_{-\alpha_r}Y_{r,2}})\Theta_{s_rs_t}\\
		&+(\tfrac{\delta_{i_r}^0\delta_{i_t}^0(Y_{t}-Y_{t,2})}{Y_rY_t^2Y_{t,2}^2}+\tfrac{\delta_{i_r}^0\delta_{i_{r,1}}^0(P_{t,2}(s_ts_r(\bi))-P_{t}(s_ts_r(\bi)))}{Y_rY_{r,1}})\Theta_{s_ts_r}
		\\
		&+(\tfrac{\delta_{i_{r}}^0\delta_{i_{t}}^0(2-Y_{t,2})(1-Y_{r,1})(Y_t-Y_{t,1})(Y_{t,1}Y_{r,2}+Y_{r,2}Y_t-Y_rY_t)}{Y_{r,1}^2Y_{r,2}^2Y_t^2Y_{t,1}^2}+\tfrac{\delta_{i_{t,2}}^0(P_{r,1}(s_r(\bi))P_t(s_r(\bi))-P_{r,2}(s_r(\bi))P_{t,1}(s_r(\bi)))}{Y_{t,2}})\Theta_r
		\\
		&+(\tfrac{\delta_{i_r}^0\delta_{i_t}^0(Y_{t,1}-1)(Y_t^2+Y_tY_{r,2}+Y_{t,2}Y_{r,1}+2Y_rY_{r,1})}{Y_rY_{r,1}Y_tY_{t,1}^2Y_{t,2}^2}+\tfrac{\delta_{i_{r,2}}^0(P_{r,1}(s_t(\bi))P_{t,2}(s_t(\bi))-P_r(s_t(\bi))P_{t,1}(s_t(\bi)))}{Y_{r,2}}\\
		&+\tfrac{\delta_{i_r}^0\delta_{i_{r,1}}^0(Y_r-Y_{r,1})P_t(s_rs_t(\bi))}{Y_r^2Y_{r,1}^2})\Theta_t]\epsilon(\bi).						\end{align*}
	There are eight cases:			
	
	$Case~1\colon i_r=i_t=0$. Using (\ref{P_r}), we have
	\begin{align*}
		&(\Psi_r\Psi_t\Psi_r\Psi_t\Psi_r\Psi_t-\Psi_t\Psi_r\Psi_t\Psi_r\Psi_t\Psi_r)\epsilon(\bi)\\
		&=[(\tfrac{Y_{t,2}-Y_{t}}{Y_t^2Y_{t,2}^2}+\tfrac{\tfrac{Y_t-1}{Y_t^2}-\tfrac{Y_{t,2}-1}{Y_{t,2}^2}}{Y_{r,1}})\Theta_{s_rs_ts_r}+(\tfrac{(Y_r-Y_{r,2})(3-3Y_{t,1}+Y_{t,1}^2)}{Y_r^2Y_{r,2}^2}+\tfrac{\tfrac{Y_{r,2}-1}{Y_{r,2}^2}-\tfrac{Y_r-1}{Y_r^2}}{Y_{t,1}})\Theta_{s_ts_rs_t}\\
		&+(\tfrac{Y_{r,2}-Y_{t,1}}{Y_{r,2}Y_t^2Y_{t,1}^2}+\tfrac{\tfrac{Y_{t,1}-1}{Y_{t,1}^2}-\tfrac{Y_t-1}{Y_t^2}}{Y_{-\alpha_r}Y_{r,2}})\Theta_{s_rs_t}+(\tfrac{Y_{t}-Y_{t,2}}{Y_rY_t^2Y_{t,2}^2}+\tfrac{\tfrac{Y_{t,2}-1}{Y_{t,2}^2}-\tfrac{Y_t-1}{Y_t^2}}{Y_rY_{r,1}})\Theta_{s_ts_r}\\
		&
		+(\tfrac{(2-Y_{t,2})(1-Y_{r,1})(Y_t-Y_{t,1})(Y_{t,1}Y_{r,2}+Y_tY_{r,2}-Y_rY_t)}{Y_{r,1}^2Y_{r,2}^2Y_t^2Y_{t,1}^2}+\tfrac{\tfrac{(Y_{r,1}-1)(Y_t-1)}{Y_{r,1}^2Y_t^2}-\tfrac{(Y_{r,2}-1)(Y_{t,1}-1)}{Y_{r,2}^2Y_{t,1}^2}}{Y_{t,2}})\Theta_r\\
		&
		+(\tfrac{(Y_{t,1}-1)(Y_t^2+Y_tY_{r,2}+Y_{t,2}Y_{r,1}+2Y_rY_{r,1})}{Y_rY_{r,1}Y_tY_{t,1}^2Y_{t,2}^2}+\tfrac{\tfrac{(Y_{r,1}-1)(Y_{t,2}-1)}{Y_{r,1}^{2}Y_{t,2}^{2}}-\tfrac{(Y_{r}-1)(Y_{t,1}-1)}{Y_r^{2}Y_{t,1}^{2}}}{Y_{r,2}}+\tfrac{(Y_r-Y_{r,1})(Y_t-1)}{Y_r^2Y_{r,1}^2Y_t^2})\Theta_t]\epsilon(\bi)\\
		&=0.
	\end{align*}
	
	$Case~2\colon i_{r,1}=i_r+3i_t=0, i_r\neq0, i_{t,1}=i_r+i_t\neq0$. It is easy to see that $i_{t,2}=i_r+2i_t=-i_t\neq0, i_{r,2}=2i_r+3i_t=i_r\neq0, e\neq 2$ and $s_rs_ts_r(\bi)_t=i_{t,3}=i_r+2i_t=-i_t, s_rs_ts_r(\bi)_{t,2}=i_t$. Thus, by (\ref{Psitheta}) and (\ref{P_r}), we have
	\begin{align*}
		&(\Psi_r\Psi_t\Psi_r\Psi_t\Psi_r\Psi_t-\Psi_t\Psi_r\Psi_t\Psi_r\Psi_t\Psi_r)\epsilon(\bi)\\
		&=\tfrac{P_t(s_rs_ts_r(\bi))-P_{t,2}(s_rs_ts_r(\bi))}{Y_{r,1}}\Theta_{s_rs_ts_r}\epsilon(\bi)\\
		&=\begin{cases}[Y_t-Y_{t,2}(Y_{t,2}-1)^{-1}]Y_{r,1}^{-1}\Psi_{r}\Psi_t\Psi_r\epsilon(\bi) &\text{if $i_r=-3, i_t=1, e\neq 2,3$};\\
			[Y_t(Y_t-1)^{-1}-Y_{t,2}]Y_{r,1}^{-1}\Psi_{r}\Psi_t\Psi_r\epsilon(\bi) &\text{if $i_r=3, i_t=-1, e\neq 2,3$};\\
			0 & \text{else}
		\end{cases}\\
		&=\begin{cases}\tfrac{1}{1-Y_{\alpha_r+2\alpha_t}}\Psi_{r}\Psi_t\Psi_r\epsilon(\bi) &\text{if $i_r=-3, i_t=1, e\neq 2,3$};\\
			\tfrac{1}{Y_t-1}\Psi_{r}\Psi_t\Psi_r\epsilon(\bi) &\text{if $i_r=3, i_t=-1, e\neq 2,3$};\\
			0 & \text{else}
		\end{cases}\\
		&=\begin{cases}(1-Y_{-\alpha_r-2\alpha_t})\Psi_{r}\Psi_t\Psi_r\epsilon(\bi) &\text{if $i_r=-3, i_t=1, e\neq 2,3$};\\
			(Y_{-\alpha_t}-1)\Psi_{r}\Psi_t\Psi_r\epsilon(\bi) &\text{if $i_r=3, i_t=-1, e\neq 2,3$};\\
			0 & \text{else}.
		\end{cases}
	\end{align*}
	
	$Case~3\colon i_{t,1}=i_r+i_t=0, i_r\neq0, i_{r,1}=i_r+3i_t\neq0$. 
	It is easy to see that $i_{t,2}=i_r+2i_t=i_t\neq 0, i_{r,2}=2i_r+3i_t=i_t\neq0, e\neq 2$ and $s_ts_rs_t(\bi)_t=i_{r,3}=2i_r+3i_t=-i_r, s_ts_rs_t(\bi)_{r,2}=i_r$. Thus, by (\ref{Psitheta}) and (\ref{P_r}), we have
	\begin{align*}
		&(\Psi_r\Psi_t\Psi_r\Psi_t\Psi_r\Psi_t-\Psi_t\Psi_r\Psi_t\Psi_r\Psi_t\Psi_r)\epsilon(\bi)\\
		&=\tfrac{P_{r,2}(s_ts_rs_t(\bi))-P_{r}(s_ts_rs_t(\bi))}{Y_{t,1}}\Theta_{s_ts_rs_t}\epsilon(\bi)\\
		&=\begin{cases}[Y_{r,2}(Y_{r,2}-1)^{-1}-Y_r]Y_{t,1}^{-1}\Psi_{t}\Psi_r\Psi_t\epsilon(\bi) &\text{if $i_r=-i_t=1, e\neq 2$};\\
			[Y_{r,2}-Y_r(Y_r-1)^{-1}]Y_{t,1}^{-1}\Psi_{t}\Psi_r\Psi_t\epsilon(\bi) &\text{if $i_r=-i_t=-1, e\neq 2$};\\
			0 & \text{else}
		\end{cases}\\
		&=\begin{cases}\tfrac{3-3Y_{\alpha_r+\alpha_t}+Y_{\alpha_r+\alpha_t}^2}{Y_{2\alpha_r+3\alpha_t}-1}\Psi_{t}\Psi_r\Psi_t\epsilon(\bi) &\text{if $i_r=-i_t=1, e\neq 2$};\\
		\tfrac{	3-3Y_{\alpha_r+\alpha_t}+Y_{\alpha_r+\alpha_t}^2}{1-Y_{r}}\Psi_{t}\Psi_r\Psi_t\epsilon(\bi) &\text{if $i_r=-i_t=-1, e\neq 2$};\\
			0 & \text{else}
		\end{cases}\\
		&=\begin{cases}(Y_{-\alpha_t}+Y_{-\alpha_r-2\alpha_t}+Y_{-2\alpha_r-3\alpha_t}-3)\Psi_{t}\Psi_r\Psi_t\epsilon(\bi) &\text{if $i_r=-i_t=1, e\neq 2$};\\
			(3-Y_t-Y_{-\alpha_r}-Y_{\alpha_r+2\alpha_t})\Psi_{t}\Psi_r\Psi_t\epsilon(\bi) &\text{if $i_r=-i_t=-1, e\neq 2$};\\
			0 & \text{else}.
		\end{cases}
	\end{align*}

	$Case~4\colon i_{t,2}=i_r+2i_t=0, i_r\neq0$. 
	It is easy to see that $i_{r,1}=i_r+3i_t=i_t\neq0, i_{t,1}=i_r+i_t=-i_t\neq0, i_{r,2}=2i_r+3i_t=-i_t\neq0, e\neq 2$ and $s_r(\bi)_t=i_{t,1}=-i_t, s_r(\bi)_{r,1}=i_{r,2}=-i_t, s_r(\bi)_{t,1}=i_t, s_r(\bi)_{r,2}=i_{r,1}=i_t$. Thus, by (\ref{Psitheta}) and (\ref{P_r}), we have
	\begin{align*}
		&(\Psi_r\Psi_t\Psi_r\Psi_t\Psi_r\Psi_t-\Psi_t\Psi_r\Psi_t\Psi_r\Psi_t\Psi_r)\epsilon(\bi)\\
		&=(P_{r,1}(s_r(\bi))P_t(s_r(\bi))-P_{r,2}(s_r(\bi))P_{t,1}(s_r(\bi)))Y_{t,2}^{-1}\Theta_{r}\epsilon(\bi)\\
		&=\begin{cases}\tfrac{Y_{r,1}Y_t-Y_{r,2}Y_{t,1}(Y_{r,2}-1)^{-1}(Y_{t,1}-1)^{-1}}{Y_{t,2}}\Psi_{r}\epsilon(\bi) &\text{if $i_r=-2, i_t=1, e\neq 2$};\\
		\tfrac{Y_{r,1}Y_{t}(Y_{r,1}-1)^{-1}(Y_{t}-1)^{-1}-Y_{r,2}Y_{t,1}}{Y_{t,2}}	\Psi_{r}\epsilon(\bi) &\text{if $i_r=2, i_t=-1, e\neq 2$};\\
			0 & \text{else}
		\end{cases}\\
		&=\begin{cases}\tfrac{(Y_{\alpha_r+2\alpha_t}-2)[Y_{\alpha_r+\alpha_t}+(Y_{2\alpha_r+3\alpha_t}-1)Y_t]}{(Y_{\alpha_r+\alpha_t}-1)(Y_{2\alpha_r+3\alpha_t}-1)}\Psi_{r}\epsilon(\bi) &\text{if $i_r=-2, i_t=1, e\neq 2$};\\
			\tfrac{(2-Y_{\alpha_r+2\alpha_t})[Y_{t}+(Y_{\alpha_r+3\alpha_t}-1)Y_{\alpha_r+\alpha_t}]}{(Y_{t}-1)(Y_{\alpha_r+3\alpha_t}-1)}\Psi_{r}\epsilon(\bi) &\text{if $i_r=2, i_t=-1, e\neq 2$};\\
			0 & \text{else}
		\end{cases}\\
		&=\begin{cases}(Y_{\alpha_r+2\alpha_t}-2)(Y_{-\alpha_r-\alpha_t}Y_{-2\alpha_r-3\alpha_t}-Y_{-\alpha_r})\Psi_{r}\epsilon(\bi) &\text{if $i_r=-2, i_t=1, e\neq 2$};\\
			(2-Y_{\alpha_r+2\alpha_t})(Y_{-\alpha_t}Y_{-\alpha_r-3\alpha_t}-Y_r)\Psi_{r}\epsilon(\bi) &\text{if $i_r=2, i_t=-1, e\neq 2$};\\
			0 & \text{else}.
		\end{cases}
	\end{align*}

	$Case~5\colon i_{r,2}=2i_r+3i_t=0, i_r\neq0$. 
	It is easy to see that $i_{r,1}=i_r+3i_t=-i_r\neq0, i_{t,1}=i_r+i_t\neq0, i_{t,2}=i_r+2i_t=-i_{t,1}\neq0$ and $s_t(\bi)_{r,1}=i_r, s_t(\bi)_{t,2}=i_{t,1}, s_t(\bi)_r=i_{r,1}=-i_r, s_t(\bi)_{t,1}=i_{t,2}=-i_{t,1}$. Thus, by (\ref{Psitheta}) and (\ref{P_r}), we have
	\begin{align*}
		&(\Psi_r\Psi_t\Psi_r\Psi_t\Psi_r\Psi_t-\Psi_t\Psi_r\Psi_t\Psi_r\Psi_t\Psi_r)\epsilon(\bi)\\
		&=\tfrac{P_{r,1}(s_t(\bi))P_{t,2}(s_t(\bi))-P_r(s_t(\bi))P_{t,1}(s_t(\bi))}{Y_{r,2}}\Theta_{t}\epsilon(\bi)\\
		&=\begin{cases}[\tfrac{Y_{r,1}^2Y_{t,2}^2}{(Y_{r,1}-1)(Y_{t,2}-1)}-\tfrac{Y_{r}^2Y_{t,1}^2}{(Y_{r}-1)(Y_{t,1}-1)}]Y_{r,2}^{-1}(\Psi_{t}+\tfrac{1}{Y_t})\epsilon(\bi) &\text{if $i_r=1,i_{t}=0, e=2$};\\
			(\tfrac{Y_{r,1}Y_{t,2}}{Y_{r,1}-1}-\tfrac{Y_{r}Y_{t,1}}{Y_{t,1}-1})Y_{r,2}^{-1}\Psi_{t}\epsilon(\bi) &\text{if $i_r=1,i_{t}=2, e=4$};\\
			(\tfrac{Y_{r,1}Y_{t,2}}{Y_{t,2}-1}-\tfrac{Y_{r}Y_{t,1}}{Y_{r}-1})Y_{r,2}^{-1}\Psi_{t}\epsilon(\bi) &\text{if $i_r=-1,i_{t}=2, e=4$};\\
			(\tfrac{Y_{r,1}}{Y_{r,1}-1}-Y_r)Y_{r,2}^{-1}\Psi_{t}\epsilon(\bi) &\text{if $i_r=1, 3i_{t}=-2, e\neq 2, 4$};\\
			(Y_{r,1}-\tfrac{Y_{r}}{Y_{r}-1})Y_{r,2}^{-1}\Psi_{t}\epsilon(\bi) &\text{if $i_r=-1,3i_{t}=2, e\neq 2, 4$};\\
			(\tfrac{Y_{t,2}}{Y_{t,2}-1}-Y_{t,1})Y_{r,2}^{-1}\Psi_{t}\epsilon(\bi) &\text{if $i_r=3, i_{t}=-2, e\neq 2, 3,4$},\\
			(Y_{t,2}-\tfrac{Y_{t,1}}{Y_{t,1}-1})Y_{r,2}^{-1}\Psi_{t}\epsilon(\bi) &\text{if $i_r=-3, i_{t}=2, e\neq 2, 3,4$};\\
			0 & \text{else}
		\end{cases}\\
		&=\begin{cases}\tfrac{(2-Y_{t})[Y_{\alpha_r+3\alpha_t}Y_{\alpha_r+2\alpha_t}+Y_rY_{\alpha_r+\alpha_t}(1-Y_t)^2]}{(1-Y_{\alpha_r+3\alpha_t})(1-Y_{\alpha_r+2\alpha_t})}(Y_t\Psi_{t}+1)\epsilon(\bi) &\text{if $i_r=1,i_{t}=0, e=2$};\\
			\tfrac{Y_{\alpha_r+\alpha_t}-Y_{\alpha_r+3\alpha_t}}{(Y_{\alpha_r+\alpha_t}-1)(Y_{\alpha_r+3\alpha_t}-1)}	\Psi_{t}\epsilon(\bi) &\text{if $i_r=1,i_{t}=2, e=4$};\\
			\tfrac{Y_{r}-Y_{\alpha_r+2\alpha_t}}{(Y_r-1)(Y_{\alpha_r+2\alpha_t}-1)}	\Psi_{t}\epsilon(\bi) &\text{if $i_r=-1,i_{t}=2, e=4$};\\
			\tfrac{1}{Y_{\alpha_r+3\alpha_t}-1}\Psi_t\epsilon(\bi)& \text{if $i_r=1,3i_{t}=-2, e\neq 2, 4$};\\
			\tfrac{1}{1-Y_{r}}\Psi_t\epsilon(\bi)& \text{if $i_r=-1, 3i_{t}=2, e\neq 2, 4$},\\
			\tfrac{1}{Y_{\alpha_r+2\alpha_t}-1}\Psi_t\epsilon(\bi) &\text{if $i_r=3,i_{t}=-2, e\neq 2,3,4$};\\
			\tfrac{1}{1-Y_{\alpha_r+\alpha_t}}\Psi_t\epsilon(\bi) &\text{if $i_r=-3,i_{t}=2, e\neq 2,3,4$};\\
			0 & \text{else}
			\end{cases}\\
			&=\begin{cases}(2-Y_{t})[Y_{-\alpha_r-2\alpha_t}(Y_{-\alpha_r-3\alpha_t}+Y_r)\\
				-Y_rY_{-2\alpha_r-3\alpha_t}](Y_t\Psi_{t}+1)\epsilon(\bi) &\text{if $i_r=1,i_{t}=0, e=2$};\\
				(Y_{-\alpha_r-\alpha_t}-Y_{-\alpha_r-3\alpha_t})	\Psi_{t}\epsilon(\bi) &\text{if $i_r=1,i_{t}=2, e=4$};\\
				(Y_{-\alpha_r}-Y_{-\alpha_r-2\alpha_t})	\Psi_{t}\epsilon(\bi) &\text{if $i_r=-1,i_{t}=2, e=4$};\\
				(Y_{-\alpha_r-3\alpha_t}-1)\Psi_t\epsilon(\bi)& \text{if $i_r=1,3i_{t}=-2, e\neq 2, 4$};\\
				(1-Y_{-\alpha_r})\Psi_t\epsilon(\bi)& \text{if $i_r=-1, 3i_{t}=2, e\neq 2, 4$};\\
				(Y_{-\alpha_r-2\alpha_t}-1)\Psi_t\epsilon(\bi) &\text{if $i_r=3,i_{t}=-2, e\neq 2,3,4$};\\
				(1-Y_{-\alpha_r-\alpha_t})\Psi_t\epsilon(\bi) &\text{if $i_r=-3,i_{t}=2, e\neq 2,3,4$};\\
				0 & \text{else}.
		\end{cases}
	\end{align*}
	
	$Case~6\colon i_{r}=0, i_{t,2}=2i_t=0, i_t\neq0$. 
	It is easy to see that $i_{r,1}=i_r+3i_t=i_t\neq0, i_{t,1}=i_r+i_t=i_t\neq0, i_{r,2}=2i_r+3i_t=i_{t}\neq0$ and $s_r(\bi)_{r,1}=i_{r,2}=i_t, s_r(\bi)_{t}=i_{t,1}=i_t, s_r(\bi)_{r,2}=i_{r,1}=-i_t, s_r(\bi)_{t,1}=i_{t}$. Thus, by (\ref{Psitheta}) and (\ref{P_r}), we obtain		
	\begin{align*}
		&(\Psi_r\Psi_t\Psi_r\Psi_t\Psi_r\Psi_t-\Psi_t\Psi_r\Psi_t\Psi_r\Psi_t\Psi_r)\epsilon(\bi)\\
		&=\tfrac{P_{r,1}(s_r(\bi))P_{t}(s_r(\bi))-P_{r,2}(s_r(\bi))P_{t,1}(s_r(\bi))}{Y_{t,2}}\Theta_{r}\epsilon(\bi)\\
		&=\begin{cases}(\tfrac{Y_{r,1}^2Y_{t,2}^2}{(Y_{r,1}-1)(Y_{t,2}-1)}-\tfrac{Y_{r}^2Y_{t,1}^2}{(Y_{r}-1)(Y_{t,1}-1)})Y_{t,2}^{-1}(\Psi_{r}+\tfrac{1}{Y_r})\epsilon(\bi) &\text{if $i_r=0, i_t=1, e=2$};\\
			0 & \text{else}
		\end{cases}\\
		&=\begin{cases}\tfrac{(Y_{\alpha_r+2\alpha_t}-2)[Y_{2\alpha_r+3\alpha_t}Y_{\alpha_r+\alpha_t}+Y_tY_{\alpha_r+3\alpha_t}(1-Y_r)]}{(1-Y_{\alpha_r+\alpha_t})(1-Y_{2\alpha_r+3\alpha_t})}(Y_r\Psi_{r}+1)\epsilon(\bi) &\text{if $i_r=0, i_{t}=1, e=2$}; \\
			0 & \text{else}
		\end{cases}\\
		&=\begin{cases}(Y_{\alpha_r+2\alpha_t}-2)[Y_{-2\alpha_r-3\alpha_t}(Y_{-\alpha_r-\alpha_t}+Y_{-\alpha_t})\\
			-Y_{-\alpha_r}Y_{-\alpha_t}](Y_r\Psi_{r}+1)\epsilon(\bi) &\text{if $i_r=0, i_{t}=1, e=2$}; \\
			0 & \text{else}.
		\end{cases}
	\end{align*}
	
	$Case~7\colon i_{r,1}=0, i_{t,1}=0, i_r\neq 0$. It is easy to see that $2i_t=0 ,i_{t,2}=i_r+2i_t=i_r\neq0, i_{r,2}=2i_r+3i_t=i_r\neq 0$ and $s_rs_ts_r(\bi)_t=i_{t,3}=i_t, s_rs_ts_r(\bi)_{t,2}=i_t, s_ts_rs_t(\bi)_{r,2}=i_r=i_t, s_ts_rs_t(\bi)_r=i_{r,3}=i_t$. Thus, by (3)-(4), (\ref{P_r}) and (\ref{Psitheta}), we obtain		
	\begin{align*}
		&(\Psi_r\Psi_t\Psi_r\Psi_t\Psi_r\Psi_t-\Psi_t\Psi_r\Psi_t\Psi_r\Psi_t\Psi_r)\epsilon(\bi)\\
		&=[\tfrac{P_{t}(s_rs_ts_r(\bi))-P_{t,2}(s_rs_ts_r(\bi))}{Y_{r,1}}\Theta_{s_rs_ts_r}+\tfrac{P_{r,2}(s_ts_rs_t(\bi))-P_r(s_ts_rs_t(\bi))}{Y_{t,1}}\Theta_{s_ts_rs_t}]\epsilon(\bi)\\
		&=[\tfrac{P_{t}(s_rs_ts_r(\bi))-P_{t,2}(s_rs_ts_r(\bi))}{Y_{r,1}}\Psi_r(\Psi_t+\tfrac{1}{Y_t})\Psi_r+\tfrac{P_{r,2}(s_ts_rs_t(\bi))-P_r(s_ts_rs_t(\bi))}{Y_{t,1}}\Psi_t(\Psi_r+\tfrac{1}{Y_r})\Psi_t]\epsilon(\bi)\\
		&=[\tfrac{P_{t}(s_rs_ts_r(\bi))-P_{t,2}(s_rs_ts_r(\bi))}{Y_{r,1}}(\Psi_{r}\Psi_t\Psi_r+\tfrac{P_r(\bi)}{Y_{t,1}})+\tfrac{P_{r,2}(s_ts_rs_t(\bi))-P_r(s_ts_rs_t(\bi))}{Y_{t,1}}(\Psi_{t}\Psi_r\Psi_t+\tfrac{P_t(\bi)}{Y_{r,1}})]\epsilon(\bi)\\
		&=\begin{cases}{\tiny[\tfrac{\tfrac{Y_t^2}{Y_t-1}-\tfrac{Y_{t,2}^2}{Y_{t,2}-1}}{Y_{r,1}}(\Psi_{r}\Psi_t\Psi_r+\tfrac{Y_r^2}{Y_{t,1}(Y_r-1)^2})
				+\tfrac{\tfrac{Y_{r,2}^2}{Y_{r,2}-1}-\tfrac{Y_r^2}{Y_r-1}}{Y_{t,1}}(\Psi_{t}\Psi_r\Psi_t+\tfrac{Y_t^2}{Y_{r,1}(Y_t-1)})]\epsilon(\bi)} &\text{if $i_r=i_t=1, e=2$};\\
			0 & \text{else}
		\end{cases}\\
		&=\begin{cases}(\tfrac{Y_{\alpha_r+\alpha_t}}{1-Y_{\alpha_r+2\alpha_t}}\Psi_{r}\Psi_t\Psi_r-\tfrac{(3-3Y_{\alpha_r+\alpha_t}+Y_{\alpha_r+\alpha_t}^2)Y_{\alpha_r+3\alpha_t}}{(1-Y_{\alpha_r+3\alpha_t})(1-Y_r)}\Psi_{t}\Psi_r\Psi_t\\
			+\tfrac{2Y_tY_{2\alpha_r+3\alpha_t}-Y_{\alpha_r+\alpha_t}Y_{\alpha_r+3\alpha_t}}{(1-Y_t)(1-Y_{2\alpha_r+3\alpha_t})})\epsilon(\bi) &\text{if $i_r=i_{t}=1, e=2$},\\
			0 & \text{else}
		\end{cases}\\
		&=\begin{cases}(Y_{-\alpha_t}-Y_{-\alpha_r-2\alpha_t})\Psi_{r}\Psi_t\Psi_r\\
			+(Y_{-\alpha_t}+Y_{-\alpha_r-2\alpha_t}+Y_{-2\alpha_r-3\alpha_t}\\
			-Y_{-\alpha_r}-Y_t-Y_{\alpha_r+2\alpha_t})\Psi_{t}\Psi_r\Psi_t\\
			+2Y_{-\alpha_t}Y_{-2\alpha_r-3\alpha_t}-Y_{-\alpha_r-\alpha_t}Y_{-\alpha_r-3\alpha_t})\epsilon(\bi) &\text{if $i_r=i_{t}=1, e=2$};\\
			0 & \text{else}.
		\end{cases}
	\end{align*}

	$Case~8\colon i_r=0, i_{r,1}=0, i_t\neq 0$. It is easy to see that $i_{t,1}=i_t\neq0, i_{t,2}=2i_t=-i_t\neq 0$ and $s_rs_ts_r(\bi)_t=i_{t,3}=-i_t, s_rs_ts_r(\bi)_{t,2}=i_{t}, s_rs_t(\bi)_t=i_{t,2}=-i_t, s_rs_t(\bi)_{t,1}=-i_t, s_ts_r(\bi)_{t,2}=i_t, s_ts_r(\bi)_t=-i_{t,1}=-i_t, s_t(\bi)_{r,1}=i_r=0, s_t(\bi)_{t,2}=i_{t,1}=i_t, s_{t}(\bi)_r=i_{r,1}=0, s_{t}(\bi)_{t,1}=i_{t,2}=-i_t, s_{t}(\bi)_t=-i_t$.	Thus, by (3)-(4) and (\ref{Psitheta}), we obtain						
	\begin{align*}
		&(\Psi_r\Psi_t\Psi_r\Psi_t\Psi_r\Psi_t-\Psi_t\Psi_r\Psi_t\Psi_r\Psi_t\Psi_r)\epsilon(\bi)\\
		&=\tiny{[\tfrac{P_{t}(s_rs_ts_r(\bi))-P_{t,2}(s_rs_ts_r(\bi))}{Y_{r,1}}\Theta_{s_rs_ts_r}+\tfrac{P_{t,1}(s_t(\bi))-P_{t}(s_rs_t(\bi))}{Y_{-\alpha_r}Y_{r,2}}\Theta_{s_rs_t}+\tfrac{P_{t,2}(s_ts_r(\bi))-P_t(s_ts_r(\bi))}{Y_rY_{r,1}}\Theta_{s_ts_r}}\\
		&\tiny{+(\tfrac{P_{r,1}(s_t(\bi))P_{t,2}(s_t(\bi))-P_{r}(s_t(\bi))P_{t,1}(s_t(\bi))}{Y_{r,2}}+\tfrac{(Y_r-Y_{r,1})P_{t}(s_rs_t(\bi))}{Y_r^2Y_{r,1}^2})\Theta_t]\epsilon(\bi)}\\
		&=\tiny{[\tfrac{P_{t}(s_rs_ts_r(\bi))-P_{t,2}(s_rs_ts_r(\bi))}{Y_{r,1}}(\Psi_r+\tfrac{1}{Y_r})\Psi_t(\Psi_r+\tfrac{1}{Y_r})+\tfrac{P_{t,1}(s_t(\bi))-P_{t}(s_rs_t(\bi))}{Y_{-\alpha_r}Y_{r,2}}(\Psi_r+\tfrac{1}{Y_r})\Psi_t}\\
		&\tiny{+\tfrac{P_{t,2}(s_ts_r(\bi))-P_t(s_ts_r(\bi))}{Y_rY_{r,1}}\Psi_t(\Psi_r+\tfrac{1}{Y_r})+(\tfrac{P_{r,1}(s_t(\bi))P_{t,2}(s_t(\bi))-P_{r}(s_t(\bi))P_{t,1}(s_t(\bi))}{Y_{r,2}}+\tfrac{(Y_r-Y_{r,1})P_{t}(s_rs_t(\bi))}{Y_r^2Y_{r,1}^2})\Psi_t]\epsilon(\bi)}\\
		&=\tiny{[\tfrac{P_{t}(s_rs_ts_r(\bi))-P_{t,2}(s_rs_ts_r(\bi))}{Y_{r,1}}(\Psi_{r}\Psi_t\Psi_r+\tfrac{1}{Y_{r,2}}\Psi_{r}\Psi_t+\parti_r(\tfrac{1}{Y_{r,1}})\Psi_t+\tfrac{1}{Y_r}\Psi_{t}\Psi_r+\tfrac{1}{Y_rY_{r,1}}\Psi_t)}\\
		&\tiny{+\tfrac{P_{t,1}(s_t(\bi))-P_{t}(s_rs_t(\bi))}{Y_{-\alpha_r}Y_{r,2}}(\Psi_{r}\Psi_t+\tfrac{1}{Y_r}\Psi_t)+\tfrac{P_{t,2}(s_ts_r(\bi))-P_t(s_ts_r(\bi))}{Y_rY_{r,1}}(\Psi_{t}\Psi_r+\tfrac{1}{Y_{r,1}}\Psi_t)}\\
		&\tiny{+(\tfrac{P_{r,1}(s_t(\bi))P_{t,2}(s_t(\bi))-P_{r}(s_t(\bi))P_{t,1}(s_t(\bi))}{Y_{r,2}}+\tfrac{(Y_r-Y_{r,1})P_{t}(s_rs_t(\bi))}{Y_r^2Y_{r,1}^2})\Psi_t]\epsilon(\bi)}\\
		&=\tiny{[\tfrac{P_{t}(s_rs_ts_r(\bi))-P_{t,2}(s_rs_ts_r(\bi))}{Y_{r,1}}(\Psi_{r}\Psi_t\Psi_r+\tfrac{1}{Y_{r,2}}\Psi_{r}\Psi_t+\tfrac{1}{Y_r}\Psi_{t}\Psi_r+\tfrac{1}{Y_rY_{r,2}}\Psi_t)}\\
		&\tiny{+\tfrac{P_{t,1}(s_t(\bi))-P_{t}(s_rs_t(\bi))}{Y_{-\alpha_r}Y_{r,2}}(\Psi_{r}\Psi_t+\tfrac{1}{Y_r}\Psi_t)+\tfrac{P_{t,2}(s_ts_r(\bi))-P_t(s_ts_r(\bi))}{Y_rY_{r,1}}(\Psi_{t}\Psi_r+\tfrac{1}{Y_{r,1}}\Psi_t)}\\
		&\tiny{+(\tfrac{P_{r,1}(s_t(\bi))P_{t,2}(s_t(\bi))-P_{r}(s_t(\bi))P_{t,1}(s_t(\bi))}{Y_{r,2}}+\tfrac{(Y_r-Y_{r,1})P_{t}(s_rs_t(\bi))}{Y_r^2Y_{r,1}^2})\Psi_t]\epsilon(\bi)}.
	\end{align*}
	There are three subcases:
	
	$Subcase~8.1\colon i_r=0, i_t=1, e=3$. By (\ref{P_r}), we obtain 
	\begin{align*}
		&(\Psi_r\Psi_t\Psi_r\Psi_t\Psi_r\Psi_t-\Psi_t\Psi_r\Psi_t\Psi_r\Psi_t\Psi_r)\epsilon(\bi)\\
		&=[\tfrac{Y_t-Y_{t,2}(Y_{t,2}-1)^{-1}}{Y_{r,1}}(\Psi_{r}\Psi_t\Psi_r+\tfrac{1}{Y_{r,2}}\Psi_{r}\Psi_t+\tfrac{1}{Y_r}\Psi_{t}\Psi_r+\tfrac{1}{Y_rY_{r,2}}\Psi_t)\\
		&+\tfrac{Y_{t,1}-Y_{t}}{Y_{-\alpha_r}Y_{r,2}}(\Psi_{r}\Psi_t+\tfrac{1}{Y_r}\Psi_t)
		+\tfrac{Y_{t,2}(Y_{t,2}-1)^{-1}-Y_t}{Y_rY_{r,1}}(\Psi_{t}\Psi_r+\tfrac{1}{Y_{r,1}}\Psi_t)\\
		&+(\tfrac{Y_{r,1}^{-2}(Y_{r,1}-1)Y_{t,2}(Y_{t,2}-1)^{-1}-Y_{r}^{-2}(Y_r-1)Y_{t,1}}{Y_{r,2}}+\tfrac{(Y_r-Y_{r,1})Y_t}{Y_r^2Y_{r,1}^2})\Psi_t]\epsilon(\bi)\\
		&=[\tfrac{1}{1-Y_{t,2}}(\Psi_{r}\Psi_t\Psi_r+\tfrac{1}{Y_{r,2}}\Psi_{r}\Psi_t+\tfrac{1}{Y_r}\Psi_{t}\Psi_r+\tfrac{1}{Y_rY_{r,2}}\Psi_t)\\
		&+\tfrac{Y_{t,1}-1}{Y_{r,2}}(\Psi_{r}\Psi_t+\tfrac{1}{Y_r}\Psi_t)
		+\tfrac{1}{Y_r(Y_{t,2}-1)}(\Psi_{t}\Psi_r+\tfrac{1}{Y_{r,1}}\Psi_t)\\
		&+(\tfrac{Y_{r,1}-Y_t}{Y_{r,1}^2Y_{r,2}}+\tfrac{(1-Y_r)Y_{t,1}}{Y_r^2Y_{r,2}}+\tfrac{(Y_r-Y_{r,1})Y_t}{Y_r^2Y_{r,1}^2})\Psi_t]\epsilon(\bi)\\
		&=\tfrac{1}{1-Y_{\alpha_r+2\alpha_t}}(\Psi_{r}\Psi_t\Psi_r+\Psi_r\Psi_t)\epsilon(\bi)\\
		&=(1-Y_{-\alpha_r-2\alpha_t})(\Psi_{r}\Psi_t\Psi_r+\Psi_r\Psi_t)\epsilon(\bi).
	\end{align*}
	
	$Subcase~8.2\colon i_r=0, i_t=-1, e=3$. By (\ref{P_r}), we obtain 
	\begin{align*}
		&(\Psi_r\Psi_t\Psi_r\Psi_t\Psi_r\Psi_t-\Psi_t\Psi_r\Psi_t\Psi_r\Psi_t\Psi_r)\epsilon(\bi)\\
		&=[\tfrac{Y_t(Y_t-1)^{-1}-Y_{t,2}}{Y_{r,1}}(\Psi_{r}\Psi_t\Psi_r+\tfrac{1}{Y_{r,2}}\Psi_{r}\Psi_t+\tfrac{1}{Y_r}\Psi_{t}\Psi_r+\tfrac{1}{Y_rY_{r,2}}\Psi_t)\\
		&+\tfrac{Y_{t,1}(Y_{t,1}-1)^{-1}-Y_t(Y_t-1)^{-1}}{Y_{-\alpha_r}Y_{r,2}}(\Psi_{r}\Psi_t+\tfrac{1}{Y_r}\Psi_t)+\tfrac{Y_{t,2}-Y_t(Y_t-1)^{-1}}{Y_rY_{r,1}}(\Psi_{t}\Psi_r+\tfrac{1}{Y_{r,1}}\Psi_t)\\
		&+(\tfrac{Y_{r,1}^2(Y_{r,1}-1)Y_{t,2}-Y_{r}^{-2}(Y_r-1)Y_{t,1}(Y_{t,1}-1)^{-1}}{Y_{r,2}}+\tfrac{(Y_r-Y_{r,1})Y_t(Y_t-1)^{-1}}{Y_r^2Y_{r,1}^2})\Psi_t]\epsilon(\bi)\\
		&=[\tfrac{1}{Y_{t}-1}(\Psi_{r}\Psi_t\Psi_r+\tfrac{1}{Y_{r,2}}\Psi_{r}\Psi_t+\tfrac{1}{Y_r}\Psi_{t}\Psi_r+\tfrac{1}{Y_rY_{r,2}}\Psi_t)\\
		&+\tfrac{1}{Y_{r,2}(1-Y_t)}(\Psi_{r}\Psi_t+\tfrac{1}{Y_r}\Psi_t)+\tfrac{1}{Y_r(1-Y_t)}(\Psi_{t}\Psi_r+\tfrac{1}{Y_{r,1}}\Psi_t)\\
		&+(\tfrac{(Y_{r,1}-1)Y_{t,2}}{Y_{r,1}^2Y_{r,2}}-\tfrac{Y_{t,1}}{Y_r^2Y_{r,2}(1-Y_t)}+\tfrac{(Y_r-Y_{r,1})Y_t}{Y_r^2Y_{r,1}^2}(Y_t-1))\Psi_t]\epsilon(\bi)\\
		&=\tfrac{1}{Y_t-1}\Psi_{r}\Psi_t\Psi_r\epsilon(\bi)\\
			&=(Y_{-\alpha_t}-1)\Psi_{r}\Psi_t\Psi_r\epsilon(\bi).
	\end{align*}
	
	$Subcase~8.3\colon i_r=0, 3i_t=0, i_t\neq \pm 1$. By (\ref{P_r}), we have
	\begin{align*}&(\Psi_r\Psi_t\Psi_r\Psi_t\Psi_r\Psi_t-\Psi_t\Psi_r\Psi_t\Psi_r\Psi_t\Psi_r)\epsilon(\bi)\\
		&=(\tfrac{Y_{r,1}^{-2}(Y_{r,1}-1)-Y_r^{-2}(Y_r-1)}{Y_{r,2}}+\tfrac{Y_r-Y_{r,1}}{Y_r^2Y_{r,1}^2})\Psi_t\epsilon(\bi)\\
		&=0.
	\end{align*}
	Then the assertion follows.			
	\end{proof}

\noindent\textit{Proof of Theorem~\ref{thm:KLR basis}.}
Relations (3)--(9) were established in Lemmas~\ref{lem:Psiepsilon}--\ref{lem:r3->t}; as noted above, relations (1) and (2) follow from the definitions. It remains to show that no additional defining relations are needed.

For each $w\in\mcW$, fix a reduced expression $w=s_{r_1}s_{r_2}\cdots s_{r_m}$ and set
\[
\Psi_w:=\Psi_{r_1}\Psi_{r_2}\cdots\Psi_{r_m}.
\]
In general, $\Psi_w$ depends on the chosen reduced expression; compare \cite[Proposition 2.5]{BKW}. From the definition of the $\Psi_r$ and the triangularity with respect to the $\Theta$-basis, we have
\[
\Psi_w\epsilon(\bi)=\left(\Theta_w+\sum_{\ell(w')<\ell(w)}f_{w'}\Theta_{w'}\right)\epsilon(\bi),
\qquad f_{w'}\in\Bbbk(Y;G).
\]
Together with the decomposition (\ref{ThetaoverlineH(E)}), this triangular relation implies, by induction on $\ell(w)$, that the elements $\Psi_w$ form a basis of $\overline{\mcH}_q(\mcC)$ over $\Bbbk(Y;G)\otimes_\Bbbk\mathcal{E}(\mcC)$. On the other hand, relations (1)--(9) reduce every word in the generators to a linear combination
\[
\sum_{w\in\mcW,\,\bi\in\mcC}\Psi_w f_{w,\bi}\epsilon(\bi),
\qquad f_{w,\bi}\in\Bbbk(Y;G).
\]
The linear independence of the $\Psi_w\epsilon(\bi)$ therefore shows that there are no further relations. Hence (1)--(9) form a complete set of defining relations.\hfill$\square$

	\begin{remark}
		By the proof of Theorem \ref{thm:KLR basis}, there is a decomposition 
		\begin{equation} \label{equa:Hq decom by Psi_w}
			\overline{\mcH}_q(\mcC)=\oplus_{\bi\in \mcC, w\in \mcW}\Bbbk(Y;G)\Psi_w\epsilon(\bi)=\oplus_{\bi\in \mcC, w\in \mcW}\Psi_w\Bbbk(Y;G)\epsilon(\bi)
		\end{equation}
	\end{remark}

	\section{Cyclotomic Hecke algebras and non-graded BKR-like isomorphisms}
	
	In this section, we fix $\Lambda=(\Lambda_i)_{i\in I}\in \mathbb{N}^I$ (we follow the convention that $\mathbb{N}=\{0,1,2,\cdots\}$) with $\sum_{i\in I}\Lambda_i<\infty$.

	\subsection{The KLR-like algebra $\mcL(\mcC)$}	
	Let $\mathcal{L}(\mcC)$ be the subalgebra of $\overline{\mcH}_q(\mcC)$ generated by $$\{Y_\alpha \ | \ \alpha\in\Sigma\}\cup \{ \Psi_r \ | \ r\in [n]\} \cup\{\epsilon(\bi) \ | \ \bi\in \mcC \}.$$
	
	The following result gives a KLR-like presentation of the subalgebra $\mcL(\mcC)$.
	\begin{proposition}
		\label{prop:complete relations of nondegenerate KLR subalg}
		$(1)$ $\mcL(\mcC)=\oplus_{\bi\in \mcC, w\in \mcW}\Bbbk[Y;G]\Psi_w\epsilon(\bi)=\oplus_{\bi\in \mcC, w\in \mcW}\Psi_w\Bbbk[Y;G]\epsilon(\bi)$.
		
		$(2)$ Relations (1)-(9) of Theorem \ref{thm:KLR basis} and 
		\begin{equation}\label{equa:Ypm} Y_\alpha Y_{\beta}=Y_\beta Y_{\alpha}, \ Y_{\alpha+\beta}=Y_\alpha+Y_\beta-Y_\alpha Y_\beta, \ Y_0=0
			\end{equation}
		for $\alpha, \beta\in \Sigma, \alpha+\beta\in \Sigma\cup\{0\}$,
		are complete for $\mcL(\mcC)$.
	\end{proposition}	
	\begin{proof}
		(1)	By the decomposition (\ref{equa:Hq decom by Psi_w}), both $\oplus_{\bi\in \mcC, w\in \mcW}\Bbbk[Y;G]\Psi_w\epsilon(\bi)$ and \\ $\oplus_{\bi\in \mcC, w\in \mcW}\Psi_w\Bbbk[Y;G]\epsilon(\bi)$ are $\Bbbk$-subspaces of $\overline{\mcH}_q(\mcC)$. Meanwhile, the relations (1)-(9) of Theorem \ref{thm:KLR basis} show that they are closed under multiplication, and thus subalgebras of $\overline{\mcH}_q(\mcC)$. By the definition of $\mcL(\mcC)$, we have $\oplus_{\bi\in \mcC, w\in \mcW}\Bbbk[Y;G]\Psi_w\epsilon(\bi)=\oplus_{\bi\in \mcC, w\in \mcW}\Psi_w\Bbbk[Y;G]\epsilon(\bi)=\mcL(\mcC)$ in $\overline{\mcH}_q(\mcC)$. 
		
		(2) Assume that $\mcL'$ is the associative unital $\Bbbk$-algebra generated by 
		$$\{Y_\alpha\ | \ \alpha\in\Sigma\}\cup \{ \Psi_r \ | \ r\in [n]\} \cup\{\epsilon(\bi) \ | \ \bi\in \mcC \}.$$
		with the relations (1)-(9) of Theorem \ref{thm:KLR basis} and (\ref{equa:Ypm}). Sending each generator of $\mcL'$ to the corresponding element of $\overline{\mcH}_q(\mcC)$ defines an algebra homomorphism. Using relations (1)--(9) of Theorem~\ref{thm:KLR basis}, every element of $\mcL'$ can be written as $\sum_{w\in\mcW, \bi\in \mcC} g_{w,\bi}\Psi_w \epsilon(\bi)$ with $g_{w,\bi}\in \Bbbk[Y;G]$. Since $\{\Psi_w \ | \ w\in \mcW\}$ is a basis of $\overline{\mcH}_q(\mcC)$ as a $\Bbbk(Y;G)\otimes_\Bbbk\mathcal{E}(\mcC)$-module, these normal forms are linearly independent after mapping to $\overline{\mcH}_q(\mcC)$. Hence the homomorphism is injective, and by (1) its image is precisely $\mcL(\mcC)$.
	\end{proof}

	\subsection{The localizations of $\mcL(\mcC)$}	
	Denote by $$\Bbbk[Y;G]^\mcW=\{f\in \Bbbk[Y;G] \ | \ w(f)=f, \ \mbox{for all} \ w\in \mcW \}.$$
	This algebra lies in the center of $\mcL(\mcC)$ by Proposition \ref{prop:complete relations of nondegenerate KLR subalg}. Let $\mathcal{K}$ be the localization of $\Bbbk[Y;G]^\mcW$ with respect to the set $\{f\in \Bbbk[Y;G]^\mcW \ | \ f\neq 0 \}$.
	
	By \cite[3.12(a)]{Lusztig89}, there is an isomorphism
	\begin{equation}\label{equa:localization1}\Bbbk[Y;G]\otimes_{\Bbbk[Y;G]^\mcW}\mathcal{K}\to \Bbbk(Y;G), \ f\otimes g\mapsto fg.\end{equation}
	Thus there is an isomorphism $\overline{\mcH}_q(\mcC)\cong \mcL(\mcC)\otimes_{\Bbbk[Y;G]^\mcW}\mathcal{K}$. 
	
	Let $\mathcal{F}$ be the localization of $\Bbbk[Y;G]^\mcW$ with respect to the set $$\{f\in \Bbbk[Y;G]^\mcW \ | \ f(0)\neq 0 \}.$$ Then one has 
	\begin{equation}\label{equa:localization2}\Bbbk[Y;G]\otimes_{\Bbbk[Y;G]^\mcW}\mathcal{F}\xrto{\cong} \mcP(Y), \ g\otimes h\mapsto gh
	\end{equation}
	where $\mcP(Y)=\Bbbk[Y;G][f^{-1}\ | \ f\in \Bbbk[Y;G], f(0)\neq 0]$ is the localization of $\Bbbk[Y;G]$ with respect to the set $\{f\in \Bbbk[Y;G] \ | \ f(0)\neq 0\}$.
	
	Let $\overline{\mcL}(\mcC)$ be the $\Bbbk$-algebra $\mcL(\mcC)\otimes_{\Bbbk[Y;G]^\mcW}\mathcal{F}$. It contains $\mcL(\mcC)$ as a $\Bbbk[Y;G]^\mcW$-subalgebra. We identify the subspace $\Bbbk[Y;G]\otimes_{\Bbbk[Y;G]^\mcW}\mathcal{F}$ of $\overline{\mcL}(\mcC)$ with $\mcP(Y)$ using (\ref{equa:localization2}). Therefore, the algebra $\overline{\mcL}(\mcC)$ can be viewed as a subalgebra of $\overline{\mcH}_q(\mcC)$ using the isomorphism $\overline{\mcH}_q(\mcC)\cong \mcL(\mcC)\otimes_{\Bbbk[Y;G]^\mcW}\mathcal{K}$.

	For $f=f(X_\alpha; \alpha\in \Sigma)\in \Bbbk[X;F] $, we write $f(q^\bi):= f(q^{i_\alpha}; \alpha\in \Sigma)$ for $\bi\in \mcC$. 
	Let $\mcP(X,\mcC)=\oplus_{\bi\in \mcC}\Bbbk[X;F][f^{-1}\ | \ f\in \Bbbk[X;F], f(q^\bi)\neq 0]\epsilon(\bi)$ be a commutative subalgebra of $ \overline{\mcH}_q(\mcC)$. Then there is an action of $\mcW$ via 
	$$w(X_\alpha)=X_{w(\alpha)}, \ w(f^{-1}\epsilon(\bi))=w(f)^{-1}\epsilon(w(\bi))$$
	where $\alpha\in \Sigma, w\in \mcW, \bi\in \mcC, f\in \Bbbk[X;F] \ \mbox{with} \ f(q^{\bi})\neq0$ (note that $w(f)(q^{w(\bi)})=f(q^\bi)$).

	\begin{lemma}\label{lem:p(X,e)=p(Y,e)} $\mcP(X,\mcC)=\oplus_{\bi\in \mcC}\mcP(Y)\epsilon(\bi)$
		in $\overline{\mcH}_q(\mcC)$. 
	\end{lemma}
	\begin{proof} The assertion follows from the identities $X_\alpha\epsilon(\bi)=(q^{i_\alpha}-q^{i_\alpha}Y_\alpha)\epsilon(\bi)$ and $Y_\alpha\epsilon(\bi)=(1-q^{-i_\alpha}X_\alpha)\epsilon(\bi)$ in $\overline{\mcH}_q(\mcC)$ for each $\alpha\in \Sigma$ and $\bi\in \mcC$.
	\end{proof}

By the defining relation (\ref{Teps}) and the decomposition (\ref{extensiondecom-Tw}), we can decompose $\overline{\mcH}_q(\mcC)$ to $\oplus_{w\in \mcW,\bi\in \mcC}\Bbbk(X;F)\epsilon(\bi)T_w$. Hence we may consider the $\Bbbk$-subspace $\mathcal{P}(X,\mcC)\otimes_\Bbbk H_q$ of $\overline{\mcH}_q(\mcC)$. This subspace is a subalgebra because it is closed under multiplication, as follows from the identities (which can be proved by (\ref{T_if}) and (\ref{Teps})):	
	\begin{align*}
	T_rX_\alpha&=X_{s_r(\alpha)}T_r+(q-1)D_r(X_\alpha), \ r\in [n], \alpha\in \Sigma,\\
	T_r f^{-1}\epsilon(\bi)&=\begin{cases}[s_r(f)^{-1}\epsilon(\bi)+(1-q)D_r(f)f^{-1}s_r(f)^{-1}\epsilon(\bi)]T_r & \text{if $i_r=0$}\\
		s_r(f^{-1}\epsilon(\bi))T_r+\tfrac{q-1}{X_r-1}(s_r(f^{-1}\epsilon(\bi)) -f^{-1}\epsilon(\bi))& \text{if $i_r\neq 0$}
	\end{cases}
\end{align*}
where $\bi\in \mcC, f\in \Bbbk[X;F]$ with $f(q^\bi)\neq0$. We will use the following consequence.

\begin{lemma}\label{lem:localization of KLR} We have $\overline{\mathcal{L}}(\mcC)\cong \mathcal{P}(X,\mcC)\otimes_\Bbbk H_q$ as subalgebras of $\overline{\mcH}_q(\mcC)$.
\end{lemma}
\begin{proof} By Lemma \ref{lem:p(X,e)=p(Y,e)}, we have $\mcP(X,\mcC)=\oplus_{\bi\in \mcC}\mcP(Y)\epsilon(\bi)$
	in $\overline{\mcH}_q(\mcC)$. By Proposition \ref{prop:complete relations of nondegenerate KLR subalg}(1) we have the decomposition $\overline{\mcL}(\mcC)=\oplus_{w\in \mcW, \bi\in \mcC}\mathcal{P}(Y)\Psi_w\epsilon(\bi)$. Note that the subalgebra $\mathcal{P}(X,\mcC)\otimes_\Bbbk H_q$ can be decomposed as $\oplus_{w\in \mcW, \bi\in \mcC}\mathcal{P}(X,\mcC)T_w\epsilon(\bi)$ using (\ref{extensiondecom-Tw}). It follows that there is an algebra isomorphism from $\overline{\mcL}(\mcC)
	$ to $ \mathcal{P}(X,\mcC)\otimes_\Bbbk H_q$ defined via 
	\begin{align*}
		Y_\alpha &\mapsto \sum_{\bi\in \mcC}(1-q^{-i_{\alpha}}X_\alpha)\epsilon(\bi),\\
		\Psi_r&\mapsto \sum\limits_{\begin{subarray}{1}\bi\in\mcC \\
				i_r\neq 0\end{subarray}}(T_r+\tfrac{q-1}{q^{i_r}-1-q^{i_r}Y_r})Q_r^{-1}(\bi)\epsilon(\bi)+\sum\limits_{\begin{subarray}{1}\bi\in\mcC \\
				i_r=0\end{subarray}}(T_r+1)Q^{-1}_r(\bi)\epsilon(\bi),\\
		\epsilon(\bi)&\mapsto \epsilon(\bi)
	\end{align*}
	for $\alpha\in \Sigma, r\in [n]$ and $\bi\in \mcC$. Recall that for $\bi\in \mcC$ and $r\in [n]$, if we write $Q_r(\bi)=\tfrac{g(Y_r)}{f(Y_r)}$ for some $f,g\in \Bbbk[Y_r]$, then $f(0)\neq0, g(0)\neq 0$. In particular, we have $\oplus_{\bi\in \mcC}\mcP(Y)T_r\epsilon(\bi)=\oplus_{\bi\in \mcC}\mcP(Y)\Psi_r\epsilon(\bi)$ for all $r\in [n]$.
\end{proof}

\begin{remark}
	We have $\oplus_{w\in \mcW,\bi\in\mcC} \mcP(Y)\Psi_w\epsilon(\bi)=\oplus_{w\in \mcW,\bi\in\mcC} \mcP(Y)T_w\epsilon(\bi)$, but in general, $\oplus_{w\in \mcW,\bi\in\mcC} \Bbbk[Y;G]\Psi_w\epsilon(\bi)\neq\oplus_{w\in \mcW,\bi\in\mcC} \Bbbk[Y;G]T_w\epsilon(\bi)$.
\end{remark}

	\subsection{The cyclotomic quotient of $\mcL(\mcC)$} We denote by 
	$$\mcL(\mcC;\Lambda):=\mcL(\mcC)/\langle Y_1^{\Lambda_{i_1}}\epsilon(\bi) \ | \ \bi\in \mcC\rangle$$
	
	As in \cite[Lemma 2.1]{Brundan-Kleshchev09}, we have the following result.
	\begin{lemma}\label{lem:deynilpotent} The elements $Y_\alpha\in \mcL(\mcC;\Lambda)$ are nilpotent for all $\alpha\in \Sigma$.
	\end{lemma}
	\begin{proof} The proof of \cite[Lemma 2.1]{Brundan-Kleshchev09} carries over. Indeed, since $\Sigma\subseteq \oplus_{r\in [n]}\mathbb{Z}\alpha_r$ and $Y_{\alpha+\beta}=Y_{\alpha}+Y_{\beta}-Y_{\alpha}Y_{\beta}, Y_0:=0$ for $\alpha, \beta\in \Sigma$ with $ \alpha+\beta\in \Sigma\cup \{0\}$, it is enough to show that $Y_r$ is nilpotent for each $r\in [n]$. For $r,t\in[n]$, let $d(r,t)$ denote the length of a shortest path connecting $r$ and $t$ in $\mfD$. We argue by induction on $d(r,1)$. If $d(r,1)=0$, then $r=1$, and the assertion follows because $Y_1^{\Lambda_{i_1}}\epsilon(\bi)=0$ in $\mcL(\mcC;\Lambda)$ for each $\bi\in \mcC$. For the induction step, assume the assertion holds whenever $d(r,1)=l$. If $d(r,1)=l+1$, then there exists $t\in [n]$ such that $d(t,1)=l$, $d(r,t)=1$ and one of the cases $r-\-t$, $\xymatrix@C10pt{r\ar@{=}[r]|{\rangle}\ar@{=}[r]&t}$ and $\xymatrix@C10pt{r\ar@3{-}[r]|{\rangle}\ar@3{-}[r]&t}$ holds. By the induction hypothesis, $Y_t$ is nilpotent.		
		
	For any $\bi\in \mcC$, if $i_r\neq 0$, by Theorem \ref{thm:KLR basis}(2)-(5) we have that	
	\begin{align*}\Psi_rY_t\Psi_r\epsilon(\bi)&=\Psi_rY_t\epsilon(s_r(\bi))\Psi_r\epsilon(\bi)\\
		&=Y_{s_r(\alpha_t)}\Psi_r^2\epsilon(\bi)\\
		&=(Y_r+Y_t-Y_rY_t)\Psi_r^2\epsilon(\bi)\\
		&=(Y_r+Y_t-Y_rY_t)Y_r^aY_{-\alpha_r}^b\epsilon(\bi)
	\end{align*}
	for some $a,b\in \{0,1\}$. Hence
	$$Y_r^{a+b+1}(Y_{-\alpha_r}-1)^b\epsilon(\bi)=\Psi_rY_t\Psi_r\epsilon(\bi)-Y_tY_r^a(1-Y_r)Y_{-\alpha_r}^b\epsilon(\bi)$$
	where we used $Y_{-\alpha_r}=Y_r(Y_{-\alpha_r}-1)$. Since $Y_t$ is nilpotent, Theorem \ref{thm:KLR basis}(4) implies that $\Psi_rY_t\Psi_r\epsilon(\bi)$ is also nilpotent. The two commuting terms on the right-hand side show that $Y_r^{a+b+1}(Y_{-\alpha_r}-1)^b\epsilon(\bi)$ is nilpotent since $$\Psi_rY_t\Psi_r\epsilon(\bi)Y_tY_r^a(1-Y_r)Y_{-\alpha_r}^b\epsilon(\bi)=Y_tY_r^a(1-Y_r)Y_{-\alpha_r}^b\epsilon(\bi)\Psi_rY_t\Psi_r\epsilon(\bi).$$ Hence $Y_r\epsilon(\bi)$ is nilpotent, because $Y_{-\alpha_r}-1$ is a unit in $\Bbbk[Y;G]$.

	Now assume that $i_r=0$. Let $Y_t^m=0$ in $\mcL(\mcC;\Lambda)$. Then by Theorem \ref{thm:KLR basis}(4), we deduce that
$$\Psi_rY_{s_r(\alpha_t)}^m\epsilon(\bi)+\partial_r(Y_t^m)\epsilon(\bi)=Y_t^m\Psi_r\epsilon(\bi)=0.
$$ Therefore, \begin{equation}\label{Y_t^m}\Psi_rY_{s_r(\alpha_t)}^mY_r\epsilon(\bi)+Y_{s_r(\alpha_t)}^m\epsilon(\bi)=Y_t^m\epsilon(\bi)=0.
\end{equation}
Multiplying on the left by $\Psi_r$ and using Theorem \ref{thm:KLR basis}(5), we have
$$(-\Psi_rY_{s_r(\alpha_t)}^{m}Y_r+\Psi_rY_{s_r(\alpha_t)}^{m})\epsilon(\bi)=0.$$
By (\ref{Y_t^m}), we have $\Psi_rY_{s_r(\alpha_t)}^{m}Y_r\epsilon(\bi)=-Y_{s_r(\alpha_t)}^{m}\epsilon(\bi)$. Therefore $$(\Psi_rY_{s_r(\alpha_t)}^{m}+Y_{s_r(\alpha_t)}^{m})\epsilon(\bi)=0.$$ Multiplying the preceding equation on the right by $Y_r$ and using (\ref{Y_t^m}) again, we obtain that $$Y_{s_r(\alpha_t)}^{m}(1-Y_r)\epsilon(\bi)=0.$$ Since $1-Y_r$ is a unit in $\mcL(\mcC;\Lambda)$, it follows that $Y_{s_r(\alpha_t)}^m\epsilon(\bi)=0$. Therefore $$Y_r\epsilon(\bi)=[Y_{s_r(\alpha_t)}+Y_t(Y_r-1)]\epsilon(\bi)$$ is nilpotent in $\mcL(\mcC;\Lambda)$.
\end{proof}

\subsection{The cyclotomic quotient of $\overline{\mcL}(\mcC)$} 	
We denote by 
$$\overline{\mcL}(\mcC;\Lambda):=\overline{\mcL}(\mcC)/\langle Y_1^{\Lambda_{i_1}}\epsilon(\bi) \ | \ \bi\in \mcC\rangle$$

\begin{lemma}
	\label{lem:L(CLambda) vs overlineL(CLambda)} There is an algebra isomorphism $\mcL(\mcC;\Lambda)\cong \overline{\mcL}(\mcC;\Lambda)$.
\end{lemma}	
	
	\begin{proof}
By Lemma~\ref{lem:deynilpotent}, the images of all $Y_\alpha$ are nilpotent in $\mcL(\mcC;\Lambda)$. Hence, if $f\in\Bbbk[Y;G]$ satisfies $f(0)\neq0$, then the image of $f$ in $\mcL(\mcC;\Lambda)$ is invertible, and its inverse is represented by an element of $\Bbbk[Y;G]$. Consequently, the canonical map
\[
\mcL(\mcC)\longrightarrow\overline{\mcL}(\mcC;\Lambda)
\]
factors through a surjective homomorphism
\[
\gamma\colon\mcL(\mcC;\Lambda)\twoheadrightarrow\overline{\mcL}(\mcC;\Lambda).
\]

Conversely, let $f\in\Bbbk[Y;G]^\mcW$ with $f(0)\neq0$. Its image under the canonical map
\[
\Bbbk[Y;G]^\mcW\hookrightarrow\mcL(\mcC)
\stackrel{\pi_1}{\twoheadrightarrow}\mcL(\mcC;\Lambda)
\]
is a unit. By the universal property of localization, there is therefore a homomorphism
\[
\zeta\colon\mathcal{F}\longrightarrow\mcL(\mcC;\Lambda).
\]
Together with $\pi_1$, this induces an algebra homomorphism
\[
\pi_1\otimes\zeta\colon
\overline{\mcL}(\mcC)
=\mcL(\mcC)\otimes_{\Bbbk[Y;G]^\mcW}\mathcal{F}
\longrightarrow\mcL(\mcC;\Lambda).
\]
For every $\bi\in\mcC$, this map sends
$Y_1^{\Lambda_{i_1}}\epsilon(\bi)$ to zero. Hence it factors through a homomorphism
\[
\xi\colon\overline{\mcL}(\mcC;\Lambda)\longrightarrow\mcL(\mcC;\Lambda).
\]
By construction, $\xi$ and $\gamma$ are inverse to each other. Thus
$\mcL(\mcC;\Lambda)\cong\overline{\mcL}(\mcC;\Lambda)$.
\end{proof}	
	
	The isomorphism in Lemma \ref{lem:localization of KLR} implies that $\overline{\mcL}(\mcC;\Lambda)$ is obtained from $\overline{\mcL}(\mcC)$ by quotienting by the ideal generated by $\prod_{i\in I}(X_1-q^{i})^{\Lambda_i}$.
	\begin{lemma} \label{prop:cndesemiration} 	$\overline{\mcL}(\mcC;\Lambda)=\overline{\mcL}(\mcC)/\langle \prod_{i\in I}(X_1-q^i)^{\Lambda_i}\rangle.$
	\end{lemma}
	\begin{proof}
Using $X_1\epsilon(\bj)=q^{j_1}(1-Y_1)\epsilon(\bj)$, we obtain
\begin{align*}
\prod_{i\in I}(X_1-q^i)^{\Lambda_i}
&=\sum_{\bj\in \mcC}\prod_{i\in I}(q^{j_1}-q^{j_1}Y_1-q^i)^{\Lambda_i}\epsilon(\bj)\\
&=\sum_{\bj\in \mcC}(-q^{j_1})^{\Lambda_{j_1}}
\prod_{\substack{i\in I\\ i\neq j_1}}
(q^{j_1}-q^{j_1}Y_1-q^i)^{\Lambda_i}
Y_1^{\Lambda_{j_1}}\epsilon(\bj).
\end{align*}
Hence
\[
\prod_{i\in I}(X_1-q^i)^{\Lambda_i}
\in \left\langle Y_1^{\Lambda_{j_1}}\epsilon(\bj)\mid \bj\in\mcC\right\rangle.
\]
Conversely, for each $\bj\in\mcC$ and $i\neq j_1$, the element
$(X_1-q^i)\epsilon(\bj)$ is invertible in $\overline{\mcL}(\mcC)$, so
\[
\prod_{\substack{i\in I\\i\neq j_1}}
[(X_1-q^i)^{\Lambda_i}]^{-1}\epsilon(\bj)
\in\overline{\mcL}(\mcC).
\]
Since
\[
Y_1\epsilon(\bj)=(1-q^{-j_1}X_1)\epsilon(\bj)
=-q^{-j_1}(X_1-q^{j_1})\epsilon(\bj),
\]
we have
\begin{align*}
Y_1^{\Lambda_{j_1}}\epsilon(\bj)
={}&(-q^{-j_1})^{\Lambda_{j_1}}
\prod_{i\in I}(X_1-q^i)^{\Lambda_i}
\prod_{\substack{i\in I\\i\neq j_1}}
[(X_1-q^i)^{\Lambda_i}]^{-1}\epsilon(\bj).
\end{align*}
Thus
\[
Y_1^{\Lambda_{j_1}}\epsilon(\bj)
\in\left\langle\prod_{i\in I}(X_1-q^i)^{\Lambda_i}\right\rangle.
\]
Therefore
\[
\left\langle Y_1^{\Lambda_{i_1}}\epsilon(\bi)\mid \bi\in\mcC\right\rangle
=
\left\langle\prod_{i\in I}(X_1-q^i)^{\Lambda_i}\right\rangle
\]
in $\overline{\mcL}(\mcC)$, and the claimed quotient description follows.
\end{proof}

	\subsection{The cyclotomic Hecke algebras} Let $\mcH_q$ be the affine Hecke algebra as defined in Section~2.
We define the {\itshape cyclotomic Hecke algebra $\mcH_q(\Lambda)$} as
$$\mcH_q(\Lambda):=\mcH_q/\langle\prod_{i\in I}(X_1-q^i)^{\Lambda_i}\rangle.$$
As in \cite[Subsection 4.1]{Brundan-Kleshchev09}, there is a system $\{e(\bi) \ | \ \bi\in I^\Sigma\}$ of mutually orthogonal idempotents in $\mcH_q(\Lambda)$ such that $1=\sum_{\bi\in I^\Sigma}e(\bi)$ and \begin{equation*}e(\bi)\mcH_q(\Lambda)=\{h\in \mcH_q(\Lambda) \ | \ (X_\alpha-q^{i_\alpha})^mh=0 \ \mbox{for all} \ \alpha\in \Sigma \ \mbox{and} \ m\gg0\}.\end{equation*}
It is easy to see that $X_\alpha e(\bi)=e(\bi)X_\alpha$ for all $\alpha\in \Sigma$ and $\bi\in I^\Sigma$.

Let $Pol^\Lambda$ denote the subalgebra of $\mcH_q(\Lambda)$ generated by (the images of) $\{X_\alpha\ | \ \alpha\in \Sigma\}$. Then $Pol^\Lambda e(\bi)$ is an algebra with identity element $e(\bi)$. If $i_\alpha\neq 0$, then the element $(1-X_\alpha) e(\bi)=(q^{i_\alpha}-X_\alpha+1-q^{i_\alpha})e(\bi)$ is a unit in $Pol^\Lambda e(\bi)$ by the nilpotency of $X_\alpha-q^{i_\alpha}$ for $\alpha \in \Sigma$. In this case, we write $(1-X_\alpha)^{-1}e(\bi)$ for its inverse in $Pol^\Lambda e(\bi)$. More generally, for a polynomial $f\in \Bbbk[X;F]$ with $f(q^{\bi})\neq0$, $fe(\bi)$ is a unit in $Pol^\Lambda e(\bi)$, in this case, we write $f^{-1}e(\bi)$ for its inverse.

\begin{lemma} \label{lem:Tre} For $r\in [n]$ and $\bi\in I^{\Sigma}$, we have
	\begin{equation*} \label{Tre-eTr} T_re(\bi)=\begin{cases}
	e(\bi)T_r & \text{if $i_r=0$},\\
	e(s_r(\bi))T_r+\tfrac{1-q}{1-X_r}(e(s_r(\bi))-e(\bi)) & \text{if $i_r\neq 0$}.
	\end{cases}
	\end{equation*}
\end{lemma}
\begin{proof}  
If $i_r=0$, then $q^{i_{s_r(\alpha)}}=q^{i_\alpha}$ for every $\alpha\in\Sigma$. Thus, by (\ref{T_if}), we have
	$$(X_\alpha-q^{i_\alpha})^mT_re(\bi)=T_r(X_{s_r(\alpha)}-q^{i_{s_r(\alpha)}})^me(\bi)+(1-q)D_r((X_\alpha-q^{i_\alpha})^m)e(\bi)=0$$
	for some positive integer $m\gg0$. Therefore $T_re(\bi)\in e(\bi)\mcH_q(\Lambda)$ and then $T_re(\bi)=e(\bi)T_re(\bi)=e(\bi)T_r.$
	
	If $i_r\neq0$, then by (\ref{T_if}), for every $\alpha\in\Sigma$,
	$$(X_\alpha-q^{i_{s_r(\alpha)}})^m[T_r(1-X_r)+1-q]e(\bi)
	=[T_r(1-X_r)+1-q](X_{s_r(\alpha)}-q^{i_{s_r(\alpha)}})^me(\bi)=0$$
	for some positive integer $m\gg0$. Therefore $[T_r(1-X_r)+(1-q)]e(\bi)\in e(s_r(\bi))\mcH_q(\Lambda)$ and then
	\begin{align*}[T_r(1-X_r)+(1-q)]e(\bi)&=e(s_r(\bi))[T_r(1-X_r)+(1-q)]e(\bi)\\
	&=e(s_r(\bi))T_r(1-X_{r})e(\bi).
	\end{align*}
	Multiplying on the right by $(1-X_r)^{-1}e(\bi)$ gives $$T_re(\bi)=e(s_r(\bi))T_re(\bi)-(1-q)
	(1-X_r)^{-1}e(\bi).$$
	Similarly,
	$$e(s_r(\bi))T_r=e(s_r(\bi))T_re(\bi)-(1-q)
	(1-X_r)^{-1}e(s_r(\bi)).$$
	Therefore $T_re(\bi)=e(s_r(\bi))T_r+\tfrac{1-q}{1-X_r}(e(s_r(\bi))-e(\bi))$, which proves the assertion.\end{proof}

Let $$e(\mcC):=\sum_{\bi\in \mcC}e(\bi)\in \mcH_q(\Lambda).$$
Then $e(\mcC)$ is a central idempotent in $\mcH_q(\Lambda)$ by Lemma~\ref{lem:Tre}. We can view $e(\mcC)$ as the identity unit of $\mcH_q(\Lambda)e(\mcC)$.

\subsection{Non-graded BKR-like isomorphisms}	

The following result is the key step in establishing the non-graded BKR-like isomorphisms.

\begin{theorem}\label{thm:Hecke vs localization KLR} There is an algebra isomorphism $\mcH_q(\Lambda)e(\mcC)\cong \overline{\mcL}(\mcC;\Lambda)$.
\end{theorem}

\begin{proof} We use Lemma~\ref{prop:cndesemiration}.
	By Lemma~\ref{lem:localization of KLR}, we can view $\mcH_q$ as a subalgebra of $\overline{\mcL}(\mcC)=\mcP(X,\mcC)\otimes_\Bbbk H_q$. Thus the canonical homomorphisms $\mcH_q\hookrightarrow \overline{\mcL}(\mcC)\twoheadrightarrow \overline{\mcL}(\mcC;\Lambda)$ induce a homomorphism
	$$\eta \colon \mcH_q(\Lambda) \to \overline{\mcL}(\mcC;\Lambda)$$
	sending the generators to the corresponding elements of $\overline{\mcL}(\mcC;\Lambda)$.
	
	We claim that $$\eta(e(\bj))=\begin{cases}\epsilon(\bj) &\text{if $\bj\in \mcC$};
		\\
		0 & \text{otherwise.}\end{cases}$$ 
	Indeed, if $\bi\in\mcC$ and $\bj\in I^\Sigma$ with $\bi\neq\bj$, then there is $\alpha\in \Sigma$ such that $j_\alpha\neq i_\alpha$. By the construction of $e(\bj)$, there exists an integer $m\gg0$ such that $(X_\alpha-q^{j_\alpha})^me(\bj)=0$. Hence $$(X_\alpha-q^{j_\alpha})^m\epsilon(\bi)\eta(e(\bj))=\epsilon(\bi)\eta((X_\alpha-q^{j_\alpha})^me(\bj))=0.$$
	The assumption $j_\alpha\neq i_\alpha$ implies that the element $(X_\alpha-q^{j_\alpha})^{-1}\epsilon(\bi)\in \overline{\mcL}(\mcC;\Lambda)$. Therefore 
	$$\epsilon(\bi)\eta(e(\bj))=(X_\alpha-q^{j_\alpha})^{-m}(X_\alpha-q^{j_\alpha})^m\epsilon(\bi)\eta(e(\bj))=(X_\alpha-q^{j_\alpha})^{-m}\epsilon(\bi)0=0.$$
	If $\bj\in I^\Sigma\setminus\mcC$, then $\eta(e(\bj))=\sum_{\bi\in \mcC}\epsilon(\bi)\eta(e(\bj))=0.$ If $\bj\in\mcC$, then
	$$\eta(e(\bj))=\sum_{\bi\in \mcC}\epsilon(\bi)\eta(e(\bj))=\epsilon(\bj)\eta(e(\bj))=\epsilon(\bj)\sum_{\bi\in I^\Sigma}\eta(e(\bi))=\epsilon(\bj)\eta(1)=\epsilon(\bj).$$
	Thus $\eta|_{\mcH_q(\Lambda)e(\mcC)}\colon \mcH_q(\Lambda)e(\mcC)\to \overline{\mcL}(\mcC;\Lambda)$ is an algebra homomorphism.

	Conversely, there is a homomorphism from $\mcP(X,\mcC)$ to $\mcH_q(\Lambda)e(\mcC)$ by sending $X_\alpha$ to $X_\alpha$ and $f^{-1}\epsilon(\bi)$($f\in \Bbbk[X;F]$ and $f(q^\bi)\neq0$) to $f^{-1}e(\bi)$ for $\alpha\in \Sigma$ and $\bi\in \mcC$. Combining this map with the canonical homomorphisms $H_q\hookrightarrow \mcH_q\twoheadrightarrow \mcH_q(\Lambda)e(\mcC)$, we obtain a linear map from $\overline{\mcL}(\mcC)=\mcP(X,\mcC)\otimes_\Bbbk H_q$ to $\mcH_q(\Lambda)e(\mcC)$ which is an algebra homomorphism by Lemma~\ref{lem:Tre}. Thus it induces a homomorphism		
	$$\rho \colon \overline{\mcL}(\mcC;\Lambda)\to \mcH_q(\Lambda)e(\mcC).$$
	The map $\rho$ is surjective, and a direct check on the generators shows that $\eta \rho$ is the identity on each generator of $\overline{\mcL}(\mcC;\Lambda)$. Hence $\rho$ is an isomorphism, with inverse $\eta|_{\mcH_q(\Lambda)e(\mcC)}$.
\end{proof}
Combining Theorem~\ref{thm:Hecke vs localization KLR} with Lemma~\ref{lem:L(CLambda) vs overlineL(CLambda)}, we obtain the following non-graded BKR-like isomorphism.
\begin{corollary}\label{cor:cyclotomic-isomorphism}
	There is an algebra isomorphism $\mcH_q(\Lambda)e(\mcC)\cong \mcL(\mcC;\Lambda)$.
\end{corollary}

\subsection{Representation-theoretic interpretation}

The isomorphism in Corollary~\ref{cor:cyclotomic-isomorphism} allows representation-theoretic questions for the central summand $\mcH_q(\Lambda)e(\mcC)$ to be expressed in terms of the KLR-like generators of $\mcL(\mcC;\Lambda)$. Since $e(\mcC)$ is a central idempotent, the category of finite-dimensional $\mcH_q(\Lambda)$-modules $M$ satisfying $e(\mcC)M=M$ is naturally identified with the category of finite-dimensional $\mcH_q(\Lambda)e(\mcC)$-modules. Hence Corollary~\ref{cor:cyclotomic-isomorphism} identifies this category with the category of finite-dimensional $\mcL(\mcC;\Lambda)$-modules. In particular, the usual representation-theoretic data of the central summand, such as simple modules, projective modules, Hom-spaces and extension groups, may be studied through the KLR-like presentation. We emphasize that this is a formal consequence of the algebra isomorphism; the point here is to describe concretely what the new generators mean on modules.

Let $M$ be a finite-dimensional $\mcH_q(\Lambda)$-module with $e(\mcC)M=M$, and let $M'$ denote the corresponding $\mcL(\mcC;\Lambda)$-module. By the construction of the idempotents and by the proof of Theorem~\ref{thm:Hecke vs localization KLR}, the generalized weight-space decomposition
\[
 M=\bigoplus_{\bi\in\mcC} e(\bi)M
\]
corresponds to the idempotent decomposition
\[
 M'=\bigoplus_{\bi\in\mcC}\epsilon(\bi)M'.
\]
More precisely, $e(\bi)$ corresponds to $\epsilon(\bi)$, and on the $\bi$-summand the element $Y_\alpha$ corresponds to
\[
 1-q^{-i_\alpha}X_\alpha.
\]
Thus the nilpotent operators $Y_\alpha$ record the generalized $X_\alpha$-eigenspace structure with eigenvalue $q^{i_\alpha}$. Moreover, Theorem~\ref{thm:KLR basis}(3) gives
\[
 \Psi_r\bigl(\epsilon(\bi)M'\bigr)\subseteq \epsilon(s_r(\bi))M',
\]
so the generators $\Psi_r$ act as transition operators between the generalized weight spaces indexed by $\bi$ and $s_r(\bi)$. The cyclotomic condition becomes
\[
 Y_1^{\Lambda_{i_1}}\epsilon(\bi)M'=0.
\]
Consequently, a module in the $\mcC$-summand can be described by the spaces $\epsilon(\bi)M'$, together with the nilpotent operators $Y_\alpha$ and the operators $\Psi_r$, subject to the relations of Theorem~\ref{thm:KLR basis} and the cyclotomic relations. In this way, the isomorphism of Corollary~\ref{cor:cyclotomic-isomorphism} translates the generalized $X$-weight-space structure of a cyclotomic Hecke module into the idempotent and nilpotent-operator language of the KLR-like presentation.

Finally, the idempotents $e(\mcC)$, as $\mcC$ ranges over the $W$-orbits in $I^\Sigma$, give central summands of $\mcH_q(\Lambda)$. The idempotent $e(\mcC)$ need not be primitive, so a fixed $\mcC$-summand may be a direct sum of several blocks. Accordingly, the construction above should be viewed as a presentation of the whole central summand rather than as a block classification. Nevertheless, it provides a concrete framework for studying finite-dimensional representations of cyclotomic Hecke algebras, and hence of affine Hecke algebras through their cyclotomic quotients, by KLR-type methods even though the algebras considered here are generally not $\mathbb Z$-graded.

\vskip10pt
{\bf Data Availability Statement:} Data sharing is not applicable to this study as no datasets were generated or analysed.
\vskip5pt
{\bf Author Contribution declaration:} Fan Kong and Zhi-Wei Li contributed equally to this work. Both authors participated in the conception and development of the research, carried out the mathematical arguments and computations, and wrote and revised the manuscript. Both authors reviewed and approved the final version of the manuscript.
\vskip5pt

 {\bf Competing Interest declaration:} There are no Competing Interests.

\end{document}